\documentclass[fontsize=12pt,a4paper,headings=normal,
twoside=false,leqno,parskip=half-,abstract=true]{scrartcl}
\usepackage[english]{babel}
\usepackage{comment}
\usepackage{graphicx}
\usepackage{subcaption}
\usepackage[utf8]{inputenc}
\usepackage{hyperref}
\hypersetup{
 pdftitle={},
 pdfauthor={},
 colorlinks=true,
 linkcolor=blue,
 citecolor=blue,
 filecolor=blue,
 urlcolor=blue}

\usepackage[pagewise]{lineno}

\usepackage[format=plain,labelfont=bf,font=small]{caption}
\usepackage{xcolor}
\usepackage[arrow, matrix, curve]{xy}
\usepackage{tikz}
\usepackage{float}

\usepackage{tikz-cd}
\usepackage{caption}
\usepackage{amsmath,amsthm}
\usepackage{amssymb}
\usepackage[normalem]{ulem}

\newtheorem{theorem}{Theorem}[section]
\newtheorem{lemma}[theorem]{Lemma}
\newtheorem{cor}[theorem]{Corollary}

\theoremstyle{definition}
\newtheorem{definition}[theorem]{Definition}

\theoremstyle{remark}
\newtheorem{remark}[theorem]{Remark}

\theoremstyle{remark}
\newtheorem{obs}[theorem]{Observation}

\numberwithin{equation}{section}

\DeclareMathAlphabet{\mathpzc}{OT1}{pzc}{m}{it}

\definecolor{refkey}{rgb}{1,0,0}
\definecolor{labelkey}{rgb}{1,0,0}

\title{Finding bifurcations in mathematical
epidemiology via reaction network methods}
\author{Nicola Vassena, \thanks
{Universität Leipzig, Germany, \texttt{nicola.vassena@uni-leipzig.de}}\\
Florin Avram, \thanks{Universit\'{e} de Pau,
 France, \texttt{florin.avram@orange.fr}}
 \\Rim Adenane,      \thanks{Universit\'e Ibn-Tofail, K\'enitra, Maroc, \texttt{rim.adenane9@gmail.com}}
 }
\date{\today}

\begin{document}

\maketitle

\begin{abstract}
  Mathematical Epidemiology (ME) shares with Chemical Reaction Network Theory (CRNT) the basic mathematical structure of its dynamical systems. Despite this central similarity, methods from CRNT have been seldom applied to solving problems in ME. We explore here the applicability of CRNT methods to find bifurcations at endemic equilibria of ME models.

We adapt three CRNT methods to the features of ME. First, we prove that essentially all ME models admit Hopf bifurcations for certain monotone choices of the interaction functions. Second, we offer a parametrization of equilibria Jacobians of ME systems where few interactions are not in mass action form. Third, for a quite general class of models, we show that periodic oscillations in closed systems imply periodic oscillations when demography is added. 

Finally, we apply such results to two families of networks: a general SIR model with a nonlinear force of infection and treatment rate and a recent SIRnS model with a gradual increase in infectiousness. We give both necessary conditions and sufficient conditions for the occurrence of bifurcations at endemic equilibria of both families.
\vspace{0.5cm}

\end{abstract}
\newpage
\tableofcontents

\section{Introduction}
 \paragraph{Motivation.} Dynamical systems are a common concern for Mathematical Epidemiology (ME), ecology, population dynamics, Chemical Reaction Network Theory (CRNT), and other related domains. All these fields struggle with computing fixed points and periodic orbits, assessing their stability, and the existence of bifurcations. As well known, dynamical systems may have a priori very complex behaviors, even in two dimensions, as illustrated by the  Hilbert $16^{th}$ open problem concerning the maximum number of cycles for a planar polynomial dynamical system.  However, the natural restriction to positivity-preserving systems with a \emph{network structure} \cite{golubitsky2023dynamics} leads often to unexpectedly simple results.
 In the last five decades, researchers succeeded in creating a whole panoply of tools that elucidates and exploits the network structure of their ODE systems. This happened notably in  CRNT, where concepts  such  as \emph{Feinberg-Horn-Jackson (FHJ) complex graph} and  \emph{deficiency} \cite{horn1972general,Fei87}, \emph{injectivity} \cite{CF05}, \emph{concordance} \cite{ShiFei12}, \emph{Species-Reaction (SR) graph} \cite{Ang07, banaji2009graph} were found to be very useful. These methods also attracted the attention of the algebraic geometry community -- see for example the dedicated chapter \emph{The Classical Theory of Chemical Reactions} in the book \cite{cox2020applications} by David A. Cox -- as well as other
 ODE researchers  \cite{soliman2013stronger,fages2015inferring,degrand2020graphical}. They have however not yet impacted other applied fields that deal with network ODEs, like ME, for example.

This raises the question of whether CRNT methods may contribute towards solving problems outside their very field. This would be quite natural, since typical problems of ME, population dynamics, ecology, etc, do possess a network structure as well.

 Motivated by such a question, our paper investigates the applicability of CRNT methods in finding bifurcation points in ME systems. We discuss three CRNT methods and adapt them to the needs of ME. We showcase their usability in two large classes of examples: a family of 3d SIRS models with nonlinear force of infection and treatment rates, which we will refer to as Capasso-Ruan-Wang (CRW) models in reference to the first authors who worked on models in such family; and a $(n+2)d$ SIRnS model where the recovered individuals are compartmentalized in $n$ sub-classes, which has been recently proposed by Adreu-Villaroig, González-Parra, and Villanueva (AGP) \cite{Carlos}. For the first case study, we unveil common features of CRW models that had not been observed before; for the second we provide the first analytical results about the existence of bifurcations for $n$ large enough.


\paragraph{A historical note.} The seminal 1920 paper \cite{Lotka20} by Alfred Lotka, which introduced the famous Lotka-Volterra system
\begin{equation}\label{eq:lotkavolterra}
\begin{cases}
\dot{x}=\alpha x -\beta xy;\\
\dot{y}=\beta xy - \delta y,
\end{cases}
\end{equation}
was concerned with chemical autocatalysis. A few years later, Vito Volterra analyzed the same system in a prey-predator context \cite{volterra26}, seemingly with no knowledge of Lotka's work. Moreover, the Lotka-Volterra system \eqref{eq:lotkavolterra} can be also seen as a simple epidemiological SI model where $x$ represents the susceptible individual with birth rate $\alpha$ and $y$ the infected individual with death rate $\delta$. The infection rate $\beta$ models the interaction between susceptible and infected.
In chemical systems, monomial reaction rates as they appear in system $\eqref{eq:lotkavolterra}$ are said to follow to \emph{law of mass-action} \cite{MA64}, particularly suited for elementary reactions in a well-mixed reactor. In a metabolic and enzymatic context, rational functions like \emph{Michaelis-Menten kinetics} \cite{MM13} are preferred, while \emph{Hill}  kinetics \cite{HIll10} model ligand-binding reactions. In population dynamics, the same nonlinearities go under the name of Holling's functional response of type I (mass action), type II (Michaelis-Menten), and type III (Hill)\cite{Holling65}. Note that the nonlinearity of the system increases by considering Michaelis-Menten or Hill kinetics rather than mass action. As a consequence, bifurcations are least likely to occur under mass action \cite{Vas,VasStad}.

 \paragraph{The unifying language of reaction networks.} The successes of CRNT motivated several researchers to propose turning \emph{mass-action systems} into a unified tool for studying applied dynamical systems. This is suggested in books like \cite{Haddad,Toth}, and some recent papers that might be associated with the unifying banner of \emph{algebraic biology} \cite{pachter2005algebraic,CDSS,macauley2020case,torres2021symbolic}.
 Unfortunately, the opposite of unification is happening as well, maybe due to excessive focus on specific examples. 

\paragraph{Some differences between ME and CRNT.} One important difference is that CRNT models have typically more dimensions than ME models. Their study renders thus indispensable concepts like the \emph{stoichiometric matrix}, whose columns are the directions in which the network evolves, and the \emph{kinetic matrix}, whose columns are the vector exponents of the appearing monomials. 

In CRNT, a central role is played by conservation laws. In most ME, there is only one: the total population (more constraints may appear though in multi-strain and meta-population models, involving populations divided into several groups).
The consequence of this is that CRNT had to develop tools that require a more delicate linear algebra treatment (for example replacing the Jacobian by the \emph{reduced Jacobian}, and its determinant by the \emph{core determinant} \cite{BaCa,CoCa}), while ME researchers typically remove the dependency just by eliminating one variable.

A third important difference is that all ME models have at least two possible equilibria (fixed points). The first, the \emph{Disease-Free-Equilibrium} (DFE),
corresponds to the elimination of all compartments involving sickness, and exists for all values of the parameters. The second one, the \emph{endemic equilibrium}, bifurcates from the DFE at a well-studied transcritical bifurcation. Furthermore, the \emph{Next Generation Matrix theorem} (NGM) or \emph{$R_0$-alternative} \cite{Diek,Van,Van08} states that such transcritical bifurcation may be expressed under very general conditions as $R_0 < 1$, where the \emph{basic reproduction number} $R_0$ is defined as the spectral radius of an operator that involves only the components transversal to the boundary. The argument gets more complicated in two-strain models (for example), where there may be at least 4 equilibria (DFE, two boundary single-strain equilibria, and the endemic coexistence of the two strains). Seemingly unknown outside the field, this law is fundamental in ME, since the analysis of boundary equilibria is a central concern. On the other hand, the CRNT literature is mostly concerned with strictly positive equilibria.

 \paragraph{Some related CRNT papers.}  One relevant paper is the recent work by Murad Banaji and Balázs Boros \cite{BaBo23}, where the authors study Hopf bifurcations in small mass-action systems. 
 The authors succeed in asserting the presence/absence of Hopf bifurcations for all (!) bimolecular networks involving three species and four reactions.  In particular, the simplest three-dimensional mass-action SIRS model is also covered by this study. They used classical dynamical systems tools like the Routh-Hurwitz methods and the local Hopf bifurcation theorem.
Attempting to go beyond this small size meets well-known computation barriers. 
The task for larger networks cannot be achieved with today's computers without exploiting further the specific network structure. In this direction, there is a connected body of work again by Banaji and co-authors on \emph{inheritance} \cite{banaji23split,banaji2023bifurcation}. In essence, they address possible enlargements of a network that preserve the capacity of having multiple equilibria or periodic oscillations or more generally bifurcations. This program aims then at finding minimal `atom' networks for several features: e.g. Hopf bifurcations \cite{BaBo23} as already mentioned. Such atom networks can consequently be enlarged ad-hoc into a bigger network under consideration. In a related spirit, but considering different nonlinearities suited in a metabolic context such as Michaelis-Menten, recent work \cite{Vas,VasStad} focuses on identifying patterns for instability or oscillations at a structural level. The tool here are the so-called \emph{Child-Selections}. We will discuss these techniques further in Section \ref{sec:3RNmethods} and adapt the results to the specific needs of ME.

\paragraph{Contents.} This paper investigates the applicability of CRNT methods in ME. As a first target, we focus on finding bifurcation points at endemic equilibria, i.e., finding endemic equilibria with non-hyperbolic Jacobian. After a preliminary Section \ref{RN} where the standard concepts are introduced, we will adapt three reaction network methods to ME in Section \ref{sec:3RNmethods}. More specifically, sub-section \ref{sec:symbolichunt} focuses on the method of  \emph{Symbolic Hunt} \cite{Vas} and shows in Theorem \ref{lem:purelyimaginary} that the structure of most ME models always implies an endemic equilibrium with purely-imaginary eigenvalues of the Jacobian, for a proper choice of monotone nonlinearities. Sub-section \ref{sec:sna} recalls the classic mass-action Jacobian parametrization method \emph{Stoichiometric Network Analysis} by Clarke \cite{ClarkeSNA}, and adapts it in Theorem \ref{thm:paramquasi} to ME systems where some reactions do not follow mass-action kinetics. Sub-section \ref{sec:inheritance}
builds on Banaji's inheritance results \cite{banaji23split, banaji2023bifurcation} and shows in Theorem \ref{thm:MEinheritance} that open ME systems with demographics admit periodic orbits if the closed system without demographics admits a Hopf bifurcation. We apply these results in two complementary families of examples. Section \ref{s:SIRSgT} considers 3d SIRS \emph{Capasso-Ruan-Wang} (CRW) systems, using methods from \ref{sec:symbolichunt} and \ref{sec:sna}. We show in Theorem \ref{thm:necessarybif} that the existence of bifurcations at endemic equilibria necessarily requires the presence of either a nonlinear treatment rate or a nonlinear infection rate. Moreover, if the infection rate follows is in a quite general rational form that includes Michaelis-Menten and Monod-Haldane kinetics, then the nonlinearity of the treatment rate is necessary: Corollary \ref{cor:necessaryhopfMM}. In sub-section \ref{sec:explicitexC} we find explicit values both for zero-eigenvalue and for Hopf bifurcation. Section \ref{sec:sirns} considers the recent $(n+2)d$ AGP model \cite{Carlos}, and shows that for $n>3$ both zero-eigenvalue and Hopf bifurcations are possible, even under mass-action kinetics: Theorem \ref{thm:sirnsn} and Theorem \ref{thm:hopfsirns}. Section \ref{s:Jin} summarizes and concludes the paper.

\section{Preliminaries on reaction networks}\label{RN}
Heuristically, an ODE system with \emph{a network structure} is a system where some symbolic terms (monomials,  fractions, etc) intervene in several equations. This gives rise to a finite set of directions called stoichiometric vectors
. 
Most of the results concerning network structured ODEs originated in CRNT; more recently, several authors proposed using the more general concept of \emph{interaction network} -- see for example Ali Al-Radhawi and co-authors who propose using \emph{biological interaction network} \cite{al2023structural}. As we mostly rely on CRNT, we will use throughout the expression `reaction network'. From an abstract standpoint, a reaction network simply consists of a set of \emph{populations} that interact via a network structure. 
   `Populations' may then stand for different epidemiological statuses of individuals (e.g. susceptible, infected, recovered), or different animal species in prey-predator models or more generally trophic networks, or chemical concentrations in metabolic and chemical reaction networks.

We proceed as neutral as possible. We define a reaction network $\mathcal{N}$ as a pair of sets $\mathcal{N}={P,R}$, where the set $P$ is the set of populations, and the set $R$ is the set of the reactions. A reaction $j$ is an ordered association between nonnegative linear combinations of populations:
\begin{equation} \label{reactionj}
 j: \quad s^{j}_1X_1+...+s^{j}_{|P|}X_{|P|} \underset{j}{\longrightarrow} \tilde{s}^{j}_1X_1+...+\tilde{s}^{j}_{|P|}X_{|P|}.
\end{equation}
Here, $X_1, ..., X_{|P|}$ indicate $|P|=card(P)$ distinct populations, e.g. susceptible, infected, recovered. Relying again on chemical wording, the nonnegative integer coefficients $s^{j},\tilde{s}^j$ are called \emph{stoichiometric coefficients}: they refer to the molecular count of the reaction. Populations on the LHS are called \emph{input}, and populations on the RHS are called \emph{output}. Reactions with no input are called \emph{birth reactions} and reactions with no output are called \emph{death reactions}.

We indicate with $x\ge 0$ the nonnegative $|P|$-vector of the densities of the populations (vector inequalities are meant coordinate-wise). The population dynamics $x(t)$ follows the system of ODEs:
\begin{equation}\label{eq:standardmaineq}
    \dot{x}=g(x):=\Gamma \mathbf{r}(x),
\end{equation}
where $\Gamma$ is the $|P|\times |R|$ \emph{stoichiometric matrix} defined as
\begin{equation}
\Gamma_{kj}=\tilde{s}_k^j-s_k^j.
\end{equation}

The rank of the stoichiometric matrix, $\operatorname{rank}(\Gamma)$, defines the \emph{rank $\rho$ of the system}. Left-kernel vectors $w$ of the stoichiometric matrix $\Gamma$ 
identify \emph{conservation laws}:
$$\frac{d}{dt}(w^Tx)=w^T\dot{x}=w^T\Gamma \mathbf{r}(x)=0.$$
Any trajectory $x(t)$ lives thus in a $\rho$-dimensional space, and the rank of the Jacobian matrix $Jac$ is at most $\rho$: 
$$\operatorname{rank}(Jac)=\operatorname{rank} (g_x)=\operatorname{rank}\bigg(\Gamma \frac{\partial \mathbf{r}(x)}{\partial x}\bigg)\le \rho.$$
In particular, $Jac$ possesses $|P|-\rho$ \emph{trivial eigenvalue zero}, for any choice of $\mathbf{r}$ and $x$.

In chemistry, values $\bar{x}$ satisfying
\begin{equation}
0=\Gamma \mathbf{r}(\bar{x})
\end{equation}
are called \emph{steady-states}, while in ecology and epidemiology the more standard expression \emph{equilibrium} is used. We warn the reader that in a chemical context \emph{equilibrium} is often used to refer to a steady state satisfying more restrictive properties such as \emph{detailed-balance}, which we are not addressing here. We recall that positive equilibria $\bar{x}>0$ in epidemiology are called \emph{endemic equilibria}, since it means that the equilibrium involves also a nonzero infected population.

The $|R|$-vector $\mathbf{r}(x)$ is the vector of the \emph{reaction rates}. 
Depending on the context, `reaction rates' goes also under the name of kinetics, interaction functions, nonlinearities, etc: see for example the excellent overviews \cite{Ang,Yu,CoCa}. Regarding the parametric choice of the reaction rates, we have already  mentioned the widely-used mass action, Michaelis--Menten, and Hill kinetics. We depict them all in the \emph{general Hill nonlinearity}:
\begin{equation}\label{eq:Hill}
r_j(x):=a_j\prod_{k=1}^{|\mathbf{P}|}\Bigg(\frac{x_k^{c^k_j}}{1+b^k_j x_k^{c^k_j}}\Bigg)^{s^j_k} \quad \text{for $a_j,c_j \in \mathbb{R}_{>0}$, $b_j \in \mathbb{R}_{\ge 0}$}, .
\end{equation}
The exponent $s^j_k$ is the stoichiometric coefficient of the population $k$ as an input of the reaction $j$, and $a_j, b_j^k, c_j^k$ are parameters. Fixing  $c^k_j=1$ for all $k,j$ defines Michaelis-Menten, and further fixing  $c^k_j=1, b^k_j=0$ defines mass action. For the nonlinearity  \eqref{eq:Hill} it is worth noting that
\begin{enumerate}
    \item $r_j(x)$ is a positive function that depends only on the populations that are the input of the reaction $j$.
    \item $r_j(x)=0$ if and only if one of the densities $x_k$ of its input populations  is zero.
    \item At any strictly positive population $x>0$,  \begin{equation}\label{def:monotonechemical}\frac{\partial r_j}{\partial x_k}>0 \text{ if and only if $X_k$ is an input of reaction $j$}.\end{equation}
\end{enumerate}
The previous three conditions can be used to define a quite general class of suited nonlinearities, which has been named in the literature \emph{weakly reversible kinetics} \cite{ShiFei12} or \emph{monotone chemical functions} \cite{Vas}. Throughout this paper, we will always be concerned with nonlinearities that satisfy conditions 1-2 above, and with only few exceptions regarding the monotonicity assumption 3. 

A model typically fixes only a parametric form of $r_j$, i.e.,
\begin{equation}
r_j(x)=r_j(x,p),
\end{equation}
e.g. \eqref{eq:Hill}. However, it does not specify the quantitative values of the parameters $p$ involved. The results, then, aim mostly at two directions:
\begin{enumerate}
\item Prove that there exists $\bar{p}$ such that a certain dynamical feature occurs.\\ (\emph{Existence results})
\item Prove that a certain dynamical feature does not occur for any parameter value $p$.\\ (\emph{Exclusion results})
\end{enumerate}
Bifurcation theory \cite{GuHo84} is a powerful tool in this context: A \emph{bifurcation} is a sudden qualitative change in the system behavior under a small change in the parameter values. In particular, bifurcations stand at the borders of parameter areas where a certain dynamical feature occurs. Identifying bifurcations is then an effective method to prove the existence of a certain dynamical behavior. For example, a \emph{saddle-node} bifurcation identifies multistationarity and \emph{Hopf} bifurcation identifies periodic oscillations. An obvious necessary condition for a bifurcation to take place at an equilibrium $\bar{x}$ is that the Jacobian matrix evaluated at $\bar{x}$ is non-hyperbolic, i.e., it possesses eigenvalues with zero-real part. A first classification of bifurcations is then possible solely based on the spectrum of the Jacobian at the bifurcating equilibrium $\bar{x}$. Saddle-node bifurcations, or any bifurcation where the number of equilibria change, necessarily require a failure of the implicit function theorem, i.e. a Jacobian with a zero-eigenvalue. In contrast, a necessary spectral condition for Hopf bifurcation is purely purely-imaginary eigenvalues of the Jacobian. For brevity of presentation, this paper is only concerned with the necessary non-hyperbolic spectral condition of a bifurcation. To prove the occurrence of a specific bifurcation, and its stability configuration, further analysis on \emph{nondegeneracy} and \emph{transversality conditions} is needed, which involves higher-order terms. For clarity, we also stress that we simply say `zero-eigenvalue' bifurcations, but of course we never mean possible trivial zero-eigenvalues due to conservation laws. Zero-eigenvalue bifurcations refer only to equilibria where the rank of the Jacobian drops and it is strictly less than $\rho$.
 
Potentially, models of interest may involve several different populations interacting by even more reactions: the number of parameters grows thus exponentially with the size of the system: this poses major challenges even for reasonably small networks. The next section presents three CRNT methods to find or exclude equilibria bifurcations of \eqref{eq:standardmaineq}. 

\section{Three CRNT methods}\label{sec:3RNmethods}

This section presents three reaction network methods, which we adapt to ME analysis.
\subsection{Symbolic hunt}\label{sec:symbolichunt}
The method of \emph{Symbolic Hunt} \cite{Vas} is concerned with instances of bifurcations for positive equilibria of the system \eqref{eq:standardmaineq}. We first provide a bird's eye view of the method. The core idea is to split the bifurcation quest into three levels.
\begin{enumerate}
\item The mere existence of a positive equilibrium of \eqref{eq:standardmaineq} necessarily requires the existence of a positive right kernel vector $\bar{\mathbf{v}}>0$ of the stoichiometric matrix $\Gamma$, i.e.,
\begin{equation}\label{eq:FluxVector}
    \Gamma\bar{\mathbf{v}}>0.
\end{equation}
We call any such positive vector $\bar{\mathbf{v}}>0$ satisfying \eqref{eq:FluxVector} an \emph{equilibrium flux vector}.
\item We then consider the Jacobian matrix $Jac$ of \eqref{eq:standardmaineq},
\begin{equation}Jac:=g_x=\Gamma \frac{\partial \mathbf{r}(x)}{\partial x},
\end{equation}
at a purely structural level, i.e., as a symbolic matrix $Jac=G(\mathbf{v}')$ with  `free' symbols $\mathbf{v}'$ that correspond to partial derivatives: A symbol $\mathbf{v}'_{jk}$  replaces $\partial r_j / \partial x_k$. If we enforce the monotonicity assumption \eqref{def:monotonechemical}, then all symbols are strictly positive at a positive equilibrium. However, we underline that here we do not refer to any evaluation at an equilibrium. We refer to the Jacobian matrix in which all the partial derivatives have been replaced by symbols as \emph{symbolic Jacobian}. Any equilibria bifurcation necessarily requires a choice of such symbols $\mathbf{v'}$ so that the symbolic Jacobian $G(\mathbf{v'})$ is non-hyperbolic, i.e., possesses eigenvalues with zero real part. In particular, such first symbolic evaluation may exclude a bifurcation behavior for a large class of systems: see for instance Theorem \ref{thm:necessaryHopfthm} below.
\item Assume now that there exists an equilibrium flux vector $\bar{\mathbf{v}}$ and a choice of symbols $\bar{\mathbf{v}}'$ such that $G(\bar{\mathbf{v}}')$ is non-hyperbolic. For any given parametric rate $r_j(x,p)$ we try to find $\bar{x}>0$ and parameter values $\bar{p}$ such that
\begin{equation}\label{eq:parameterconstraints}
\begin{cases}
r_j(\bar{x},\bar{p})=\bar{v}_j\\
\partial r_j(\bar{x},\bar{p})/\partial x_k=\bar{v}'_{jk}.
\end{cases}
\end{equation}
Towards this purpose, the definition of \emph{parameter-rich kinetics} has been introduced \cite{VasStad}. In the monotone case with \eqref{def:monotonechemical} enforced, the definition characterizes the existence of positive solutions to \eqref{eq:parameterconstraints} for any choice of $\bar{x}>0$. In particular, if the network is endowed with parameter-rich kinetics, then any Symbolic Jacobian can be realized as the Jacobian of a positive equilibrium of the system, i.e., \eqref{eq:parameterconstraints} has always a solution. The class of monotone parameter-rich kinetics contains Michaelis--Menten and Hill, but it excludes mass action. For systems that include reaction rates that are in mass action form, we can still use the approach both to exclude bifurcations or as a first guidance to find a bifurcation. The solution of \eqref{eq:parameterconstraints} must be however pursued case-by-case.
\end{enumerate}

\paragraph{Symbolic Jacobians from ME admit purely imaginary eigenvalues.} We now turn to ME and we prove that a very common network pattern implies purely imaginary eigenvalues of the symbolic Jacobian. To be precise, we consider an epidemiological model to contain at least two populations, called Susceptible $S$ and Infected $I$. Furthermore, the reaction for an infection to occur requires one Susceptible and one Infected to meet. The meeting results in two Infected persons. This reaction has stoichiometric form \eqref{reactionj}
\begin{equation}\label{eq:si2i}
    S+I\quad\underset{1}{\rightarrow}\quad 2I.
\end{equation}
Assuming that the reaction rates are monotone chemical, \eqref{eq:si2i} reads into the equations as
\begin{equation*}
\begin{cases}
\dot{s}=-r_1(s,i)+ ...\\
i'=r_1(s,i)+...
\end{cases}
\end{equation*}
It is also reasonable to assume that infected people recover after a while, developing protection against the disease. In case the disease is deadly, Infected people may as well die. In both cases, we can consider a reaction of the type
\begin{equation}\label{eq:ideath}
    I\quad\underset{2}{\rightarrow}\quad ...\quad,
\end{equation}
where the dots $...$ indicate any linear combination of the populations without Infected and Susceptible. The ODEs induced by \eqref{eq:si2i} and \eqref{eq:ideath} are
\begin{equation}
\begin{cases}
\dot{s}=-r_1(s,i)+ ...\\
i'=r_1(s,i)-r_2(i)+...\\
...
\end{cases}\end{equation}

The functions $r_1(s,i)$ and $r_2(i)$ are assumed monotone increasing in $s$, $i$. Our result is:
\begin{theorem} \label{lem:purelyimaginary}
The symbolic Jacobian of any epidemiological system including reactions $1$ and $2$ admits purely imaginary eigenvalues.
\end{theorem}

\proof
Consider a network with populations $(S,I,...)$, and reactions $(1,2,...)$. Let
\begin{equation}
G(\varepsilon)=\begin{pmatrix}
    \begin{matrix}
        -v'_{1S}+\varepsilon & -v'_{1I}+\varepsilon\\
        v'_{1S}+\varepsilon & v'_{1I}-v'_{2I}+\varepsilon
    \end{matrix} &  ... \\
   ... & ... \\
\end{pmatrix},
\end{equation}
denote the symbolic Jacobian where all the symbols $v'_{jk}$ associated to any partial derivative $\partial r_j/\partial x_k$ with $j \neq 1,2$ have been rescaled as $v'_{jk}=\varepsilon$. At $\varepsilon=0$, the symbolic Jacobian $G(0)$ reads
\begin{equation}
G(0)=\begin{pmatrix}
    \begin{matrix}
        -v'_{1S} & -v'_{1I}\\
        v'_{1S} & v'_{1I}-v'_{2I}
    \end{matrix} &  \text{\Large $0_1$}  \\
   ... & \text{\Large $0_2$} \\
\end{pmatrix},
\end{equation}
where $0_1$ and $0_2$ are  zero matrices. In particular the eigenvalues of $G(0)$ are $0$ with algebraic multiplicity $|\mathbf{P}|-2$  plus the eigenvalues of the block
\begin{equation}
\begin{pmatrix}
        -v'_{1S} & -v'_{1I}\\
        v'_{1S} & v'_{1I}-v'_{2I}
    \end{pmatrix},
\end{equation}
which has determinant $v'_{1S}v'_{2I}>0$ and trace $-v'_{1S} + v'_{1I} -v'_{2I}$.

Therefore, the matrix $G(0)$ has a purely imaginary crossing of the eigenvalues, with a strictly positive distance from the
origin in the complex plane,
     whenever the three free variables are chosen to satisfy
    $$v'_{2I}=v'_{1I}-v'_{1S}.$$
  This implies,  by continuity of the eigenvalues, that the matrix $G(\varepsilon)$
  also admits a purely imaginary crossing for small enough,  positive values of the symbols $v'_{jk}$, with $j \neq 1,2$.
\endproof

\begin{remark} For systems that are parameter-rich and admit a positive endemic equilibrium, Theorem \ref{lem:purelyimaginary} implies the existence of an equilibrium with purely imaginary eigenvalues. For systems that include many reactions in mass-action form, in contrast, Theorem \ref{lem:purelyimaginary} is non-conclusive, and the existence of periodic solutions must be checked case-by-case.
\end{remark}



\subsection{Jacobian parametrization}\label{sec:sna}

The symbolic approach described in Section \ref{sec:symbolichunt} is particularly effective either to exclude bifurcations or to prove the existence of bifurcations when the reaction rates \emph{all} are sufficiently rich in parameters as - for example - Michaelis-Menten. On the other hand, ME models often assume mass-action kinetics for most reactions, especially when such reactions do not involve infected individuals.
This section reviews a CRNT parametrization method to study equilibria bifurcation of mass-action systems and adapts it to include a few reaction rates of different forms.

\paragraph{A classic CRNT method.} One standard reference is \emph{Stoichiometric Network Analysis} by B. L. Clarke \cite{ClarkeSNA}. We give the elementary basic ideas from scratch. The central observation is that for an univariate polynomial $
p(x)=ax^m$,
the derivative can be expressed as
\begin{equation}\label{eq:mader}
p'(x)=m\;ax^{m-1}=m\;ax^{m-1}\frac{x}{x}=m\frac{p(x)}{x}.
\end{equation}
This simple observation has important consequences for mass action systems
\begin{equation}\label{eq:mainma}\dot{x}=\Gamma\mathbf{r}(x),
\end{equation}
where now $\mathbf{r}(x)$ are in mass-action form
\begin{equation}
r_j(x)=a_j\prod_{k=1}^{|\mathbf{P}|}x^{s^j_k}.
\end{equation}
In fact, if $\bar{x}>0$ is a positive equilibrium of \eqref{eq:mainma}, then necessarily the reaction rates evaluated at $\bar{x}$ constitute an equilibrium flux vector $\bar{\mathbf{v}}$:
$$ 0<r(\bar{x})=\bar{\mathbf{v}} \quad \text{with}\quad \Gamma \bar{\mathbf{v}}=0.$$
Crucially, $\bar{\mathbf{v}}$ only depends on the stoichiometric matrix $\Gamma$, and not on $\bar{x}$. Thus the Jacobian of any positive equilibrium $\bar{x}>0$ can be parametrized as
\begin{equation}\label{eq:jacpar}
    \begin{split}
Jac(\bar{x})&=\Gamma\frac{ \partial \mathbf{r}(x)}{\partial x}\bigg|_{x=\bar{x}}\\
    &=\Gamma\frac{ \partial \mathbf{r}(x)}{\partial x}\bigg|_{x=\bar{x}} \operatorname{diag}(\bar{x}_k/\bar{x}_k)
\\
&=\Gamma B(\bar{\mathbf{v}})\operatorname{diag}(1/\bar{x}_k)
    \end{split}
\end{equation}
where $B(\bar{\mathbf{v}}):=\partial \mathbf{r}(x)/\partial x|_{x=\bar{x}} \operatorname{diag}(\bar{x}_k)$ just generalizes \eqref{eq:mader} and does not explicitly depend on $\bar{x}$ anymore. With a straightforward calculation, it is possible to show that $B(\bar{\mathbf{v}})$ reads
\begin{equation}
B(\bar{\mathbf{v}})=\operatorname{diag}(\bar{\mathbf{v}})K^T,
\end{equation}
where $K$ is the $|P|\times|R|$ \emph{kinetic matrix} defined by the input stoichiometric coefficients:
$$K_{kj}:=s^j_k.$$

The values $1/\bar{x}_k$ can be thought as parameters themselves, and thus $\operatorname{diag}(1/\bar{x}_k)$ can be seen as a positive diagonal parametric matrix $D$. In particular, any parametrization of the \emph{flux cone}
$$\mathcal{F}_\Gamma:=\{\mathbf{v}\in \mathbb{R}^{|\mathbf{R}|}\;\text{ s.t. } \;\Gamma \mathbf{v}=0\},$$
leads to a parametrization of all equilibrium Jacobians. The flux-cone parametrization originally used by Clarke was \emph{convex coordinates}, but ME systems often have a very natural parametrization of the flux cone, directly available at first inspection.

We finally underline that a Jacobian in the form $Jac=AD$ as \eqref{eq:jacpar}, where $A:=\Gamma B(\bar{\mathbf{v}})$ and $D:=\operatorname{diag}(1/\bar{x}_k)$, is particularly convenient to find Hopf bifurcations, even without any Routh-Hurwitz computation. In fact, a change of stability of $AD$ by tuning the parameters in $D$ leads univocally to periodic orbits, as
$$\operatorname{rank}AD=\operatorname{rank}A$$
excludes zero-eigenvalue bifurcations. See \cite{Vas24noHurw} for more details, as well as sufficient conditions and connections to the linear algebra concepts of $D$-stability and $P$-matrices.

\paragraph{An adaptation for ME.} Obviously, the above parametrization method can be identically used also in any ME system where reactions are only in mass action form. However, we generalize it for a larger class of ME systems defined as follows.
\begin{definition}[Quasi mass-action system]
An ME system is said to be \emph{quasi mass-action} if there is a single special population class $I$, called here \emph{infected},
such that the following two conditions hold true:  \begin{enumerate}
\item All reaction rates are in the class of general Hill nonlinearities \eqref{eq:Hill};
    \item All reaction rates are in mass-action form with respect to any variable $x_k\neq i$, i.e., which is not the infected $i$.  
\end{enumerate}
\end{definition}
The convenience of the Hill nonlinearity is that a generalization of \eqref{eq:mader} holds,  since
$$r_j(i)=a_j\frac{i^{c_j}}{1+{b_j}i^{c_j}}$$
implies
\begin{equation}\label{eq:HMHder}
r_j'(i) = a_j\frac{c_ji^{c_j - 1}}{(1 + b_j i^{c_j})^2}=a_j\frac{c_ji^{c_j-1}i}{(1 + b_j i^{c_j})^2i}=\frac{1}{1+b_ji^{c_j}}
\frac{c_jr_j(i)}{i}\;\in (0, c_jr_j(i)],
\end{equation}
from which we derive the following parametrization theorem.
\begin{theorem}[Parametrization of Jacobians of ME quasi mass-action systems]\label{thm:paramquasi} Consider an ME quasi mass-action system, \begin{equation}\label{eq:mainqma}
    \dot{x}=\Gamma \mathbf{r}(x),
\end{equation} and let $\bar{\mathbf{v}}>0$ be a positive equilibrium flux vector, i.e.,
$$\Gamma \bar{\mathbf{v}} = 0.$$
Define a parametric vector $\mathbf{q}\in \mathbb{R}^{|\mathbf{E}|}$  as follows:
$$q_j=\begin{cases}
q_j\ge1/c_j\quad\text{ if the reaction $j$ is not in mass-action form;}\\
1\quad\quad\;\;\quad\quad\text{ otherwise.}
    \end{cases}$$
Let $\bar{\mathbf{u}}(\bar{\mathbf{v}},\mathbf{q})\in \mathbb{R}^{|\mathbf{E}|}$ be defined component-wise as $\bar{u}_j:=\bar{v}_j/q_j.$ Then the equilibrium Jacobian $Jac(\bar{x})$ of the quasi mass-action system \eqref{eq:mainqma} can be parametrized as follows:
\begin{equation}\label{eq:parjacquasi}
Jac(\bar{x})=A(\mathbf{\bar{u}})\operatorname{diag}(1/\bar{x}_k),
\end{equation}
where $A(\mathbf{\bar{u}})$ is the matrix formed by replacing the $i^{th}$ column of $A(\mathbf{\bar{v}})=\Gamma \operatorname{diag}(\mathbf{\bar{v}})K^T$ by the column vector $\Gamma\operatorname{diag}(\bar{\mathbf{u}})K_i^T$. $K_i$ indicates here the row of the kinetic matrix $K$ corresponding to the infected population $I$.
\end{theorem}

\begin{remark}
Hill nonlinearities have more parameters than mass action, and consequently the parameterization of the Jacobian \eqref{eq:parjacquasi} involves further parameters $\mathbf{q}$. We have opted to express such parameters in terms of the Hill parameters $c_j$ to underline the relevant Michaelis--Menten case with $c_j=1$, and thus $q_j\ge 1$. Theorem \ref{thm:paramquasi} clearly suggests that general Hill systems with $c_j\ge1$ inherit bifurcations from mass action systems, via a simple perturbation argument $q_j=1+\varepsilon$. We are however not pursuing the details here.
\end{remark}

\proof[Proof of Theorem \ref{thm:paramquasi}]
For Jacobian of \eqref{eq:mainqma} we again operate as before:
\begin{equation}\label{eq:jacproofquasi}
Jac(\bar{x})=\Gamma \bigg(\frac{\partial \mathbf{r}(x)}{\partial x}\bigg|_{x=\bar{x}} \operatorname{diag}(\bar{x}_k)\bigg) \operatorname{diag}(1/\bar{x}_k).
\end{equation}

By definition of ME quasi mass-action systems, any non-mass-action form in the variable $x_k$ necessarily requires $x_k=i$. That is, the differences with the pure mass-action parametrization \eqref{eq:jacpar} are located exclusively in the $i^{th}$ column of the Jacobian.

Again as in \eqref{eq:HMHder}, for $$r_j(x)=a_j \bigg(\frac{i^c_j}{1+b_ji^{c_j}}\bigg)^{s^j_i}$$ in Hill-form in $i$ we have that
$$\frac{\partial r_j(\bar{x})}{\partial i}=s^j_i\frac{c_jr(\bar{x})}{(1+b_ji^{c_j})i}.$$

Noting that the kinetic matrix is defined precisely as $K_{kj}:=s^j_k$, and introducing the parameter-vector $\mathbf{q}$ defined by
\begin{equation}\label{eq:qmap}
q_j(b_j,c_j, \bar{i}):=\frac{(1+b_j \bar{i}^{c_j})}{c_j},
\end{equation}
we get that the $i^{th}$ column of the Jacobian \eqref{eq:jacproofquasi} reads
$$\Gamma \operatorname{diag}\bigg(\frac{\bar{v}_j}{q_j}\bigg)K_i^T \frac{1}{\bar{i}}.$$
Finally note that \eqref{eq:qmap} defines a surjective map onto $[1/{c_j},+\infty)$, for free choices of $b_j\ge 0, c_j>0$. Of course, mass action, $b_j=0, c_j=1$ fixes uniquely  $q_j=1$.

\endproof

\subsection{Inheritance}\label{sec:inheritance}
In ME, it is often possible for a single individual to `travel' through all the epidemiological classes of the model. We indicate such `travel' through epidemiological classes as $\rightsquigarrow$ to underline the difference with the reaction arrow $\rightarrow$ which takes further care of the stoichiometric coefficients. For instance, all 2d $SI$ models considers
\begin{equation}
S\quad\rightsquigarrow \quad I,
\end{equation}
or also any 3d $SIR$ models  contains
$$S\quad\rightsquigarrow \quad I \quad \rightsquigarrow \quad R$$
More in general, for abstract models with $n$ populations,
$$X_1\quad\rightsquigarrow \quad X_2 \quad \rightsquigarrow \quad ... \quad \rightsquigarrow \quad X_n.$$
Inspired by such construction, we define \emph{travel-through networks} as follows:

\begin{definition}\label{def:cyclic}
Consider a network $\mathcal{N}=\{\mathbf{P},\mathbf{R}\}$ with exactly $n$ populations $X_1,...,X_N$. We call $\mathcal{N}$ \emph{travel-through} if, possibly among other reactions, there are $n-1$ reactions $j=1,...,n-1$, whose associated \emph{stoichiometric vector} $\Gamma^j$, defined as $\Gamma^j_k:=\tilde{s}^j_k-s^j_k$, reads
\begin{equation}\Gamma^j=\begin{pmatrix}
    0\\
    \vdots\\
    -1\\
    1\\
    \vdots\\
    0\\
\end{pmatrix}\begin{tabular}{c}
\\
\\
j  \\
j+1\\
 \\
 \\
\end{tabular},
\end{equation}
where the dots are filled with zero.
\end{definition}

Based on Banaji's \emph{inheritance} results \cite{banaji23split, banaji2023bifurcation}, we establish an inheritance result for travel-through networks in ME, when demographics is added. For simplicity, we focus on mass action. Let
\begin{equation}\label{eq:MEsystem}
\dot{x}=g(x)=\Gamma \mathbf{r}(x)
\end{equation}
be any ME system in mass-action form. Let $B\neq 0$ be a nonnegative constant \emph{birth} vector $B\in \mathbf{R}^n_{\ge0}$, and let $D$ be a nonnegative diagonal \emph{death} matrix, $D\in \mathbb{R}^{n\times n}_{\ge 0}$ satisfying
\begin{equation}\label{eq:demoimplies}
    B_i \neq 0 \implies D_{ii} \neq 0.
\end{equation}  We call the system \begin{equation}\label{eq:MEdemo}
  \dot{x}=B+g(x)-Dx
\end{equation}
the associated system to \eqref{eq:MEsystem} with \emph{added demographics}.
In particular, \eqref{eq:demoimplies} requires that any population to which a birth reaction is added has also an added death reaction. We do not - however - require the converse. 

Before we state and prove out result, we recall some standard concepts \cite{GuHo84}. We call \emph{periodic orbit} any \emph{nonstationary}, i.e. non-equilibrium, trajectory $x(t)$ such that $x(t+T)=x(t)$ for some \emph{period} $T\in{\mathbb{R}}_{>0}$. 
We say that a periodic orbit is \emph{nondegenerate} if its nontrivial \emph{Floquet multipliers}, i.e. the eigenvalues of its associated \emph{Poincaré map}, do not lie on the unit circle. In analogy, a Hopf bifurcation is nondegenerate if its first \emph{Lyapunov coefficient} or \emph{focal value} is nonzero. In particular, a nondegenerate Hopf bifurcation implies the existence of a nondegenerate periodic orbit. We prove the following theorem.

\begin{theorem}\label{thm:MEinheritance}
Consider a travel-through system \eqref{eq:MEsystem} in mass-action form, which admits a nondegenerate Hopf bifurcation.
Then its associated system with added demographics \eqref{eq:MEdemo}
admits nondegenerate periodic orbits.
\end{theorem}

As already mentioned, Banaji and co-authors established a list of \emph{inheritance results} \cite{banaji23split, banaji2023bifurcation}, i.e., enlargements of a network that preserve the ability to multistationarity, oscillations, or in general local bifurcations. For the purposes of this paper, we only highlight one of such enlargements on which we build, and we refer to the original references for further reading on the topic.

\begin{theorem}[Theorem 3.2, \cite{banaji2023bifurcation} and Theorem 1, \cite{banaji23split}]\label{thm:banaji}
Let $\mathcal{N}=\{\textbf{P},\textbf{R}\}$ be a network whose associated mass-action system \eqref{eq:standardmaineq} admits a nondegenerate Hopf bifurcation (respectively a nondegenerate periodic orbit). Define the enlarged reaction set $\tilde{\textbf{R}}:=\textbf{R}\cup j$, obtained by $\textbf{R}$ with the addition of a single reaction $j$ such that
$$\operatorname{rank}(\Gamma)=\operatorname{rank}(\tilde{\Gamma}).$$
Consider the enlarged network $\tilde{\mathcal{N}}:=\{\textbf{P},\tilde{\textbf{R}}\}$. Then the associated mass-action system \eqref{eq:standardmaineq} admits a nondegenerate Hopf bifurcation (respectively a nondegenerate periodic orbit).
\end{theorem}

We now proceed with the proof of Theorem \ref{thm:MEinheritance}.
\proof[Proof of Theorem \ref{thm:MEinheritance}] Since system \eqref{eq:MEsystem} undergoes a Hopf bifurcation, there exists in particular a bifurcation parameter $\beta$ and a closed interval $[\beta_1,\beta_2]$ such that
\begin{equation}\label{eq:globalhopfconditions}
    \begin{cases}
    (i) \;\;\quad g(\bar{x}(\beta),\beta)=0\hspace{5cm}\quad\quad\quad \\
    (ii)\; \quad \operatorname{det}g_x(\bar{x}(\beta),\beta)\neq 0\hspace{5cm}\quad\text{ for $\beta\in [\beta_1,\beta_2]$}\\
    (iii) \quad \operatorname{spectrum}g_x(\bar{x}(\beta_1),\beta_1)\neq \operatorname{spectrum}g_x (\bar{x}(\beta_2),\beta_2).
    \end{cases}
\end{equation}

That is, (i) $\bar{x}(\beta)$ is a local 1-dimensional continuum of equilibria (ii) with an invertible Jacobian, and (iii) a net change of stability from end to end.

Consider at first the diagonal parametric $\tilde{D}$ defined as follows
\begin{equation}
\begin{cases}
\tilde{D}_{ii}:=d_{ii}\ge 0 \quad \Leftrightarrow \quad B_i>0\\
\tilde{D}_{ii}=0,\quad \text{otherwise},
\end{cases}
\end{equation}
and fix $B$ to be
$$0\neq B_i(\beta) = d_{ii} \bar{x}(\beta).$$
For such choice of $B$ we get that, for any choice of $\tilde{D}$,
$$(i) \quad f(\bar{x}(\beta),\beta)=B(\beta)+g(\bar{x}(\beta),\beta)-\tilde{D}\bar{x}(\beta)=0,$$
that is, $\bar{x}(\beta)$ is also a local continuum of equilibria of \eqref{eq:MEdemo}. Moreover, for $\tilde{D}_{ii}=d_{ii}=\varepsilon$ small enough, the continuity of eigenvalues with respect to the entries implies that
\begin{equation}
\begin{cases}
    (ii)\; \quad \operatorname{det}f_x(\bar{x}(\beta),\beta)\neq 0\hspace{5cm}\quad\text{ for $\beta\in [\beta_1,\beta_2]$}\\
    (iii) \quad \operatorname{spectrum}f_x(\bar{x}(\beta_1),\beta_1)\neq \operatorname{spectrum}f_x (\bar{x}(\beta_2),\beta_2).
    \end{cases}
\end{equation}
Thus we have again (i), (ii), (iii) as in \eqref{eq:globalhopfconditions} that - via intermediate value theorem - imply purely imaginary of the Jacobians. Moreover, since \eqref{eq:MEdemo} is analytic and so is the parametrization $\bar{x}(\beta)$, Theorem 4.7 of Fiedler \cite{fiedler1985index} yields periodic orbits via global Hopf bifurcation. Note that the nondegeneracy of the periodic orbit for the closed system is inherited by the periodic orbit in the open system, by continuity. To conclude for general $D$ with $$B_i\neq 0 \Rightarrow D_{ii} \neq 0,$$ note that the rank of any travel-through system is at least $n-1$, since the stoichiometric matrix contains - up to relabeling - the following $n\times (n-1)$ full-rank matrix $C$:
$$C=\begin{pmatrix}
    -1 & 0 & ... & 0\\
    1 & -1 & ...& 0\\
    \vdots & \vdots& \vdots& \vdots\\
    0& 0 & ...& -1\\
    0 & 0 & ... &1
\end{pmatrix},$$
which clearly satisfies $\operatorname{rank}C=n-1.$ Then, any travel-through system with added any birth or death from at least one population has full rank $n$: Consider for example the $n \times n $ $\tilde{C}$:
$$\tilde{C}=\begin{pmatrix}
    -1 & 0 & ... & 0 & 1\\
    1 & -1 & ...& 0 & 0\\
    \vdots & \vdots& \vdots& \vdots & \vdots\\
    0& 0 & ...& -1 & 0\\
    0 & 0 & ... &1 &0
\end{pmatrix}.$$

Note that the vector $B$ is nonzero, so there is at least one population with added both birth and death reactions. The above argument for $\tilde{D}$ guarantees that we have then a full-rank system with nondegenerate periodic orbits. Banaji's Theorem \ref{thm:banaji} further guarantees that we can add any reaction involving the existent populations to the system with $\tilde{D}$ and the resulting enlarged system still admits periodic orbits.
\endproof

\section{Application I: Capasso-Ruan-Wang SIRS epidemic models}\label{s:SIRSgT}

This section showcases our results from sections \ref{sec:symbolichunt} and \ref{sec:sna}. We apply them to conclude stability and bifurcation results for \emph{Capasso-Ruan-Wang} nonlinear SIRS epidemic models
\begin{equation}\label{Vys}
\begin{cases}
\dot{s}= \lambda  -\gamma_ss    - \beta s iF(i) 
+ \gamma_r r -\mu_s s,
\\
i'=
\beta s i F(i) -\gamma_i i -T(i)-\mu_i i,  \\
\dot{r}= \gamma_s s + \gamma_i   i +T(i) -\gamma_r r-\mu_r r  ,
\end{cases}
\end{equation}
with nonnegative, nondecreasing \emph{treatment rate} $T(i)\ge0$, $T'(i)\ge0$, and positive \emph{inhibition function} $F>0$ with $F(0)=1$. We assume positive $\beta,\gamma_i>0$ and nonnegative $\lambda,\gamma_s,\mu_s, \mu_i, \gamma_r, \mu_r \ge 0$.

 We see this model as \emph{semi-parametric} instead of symbolic because we assume specific types
of functional dependencies. Let us explain the system from a reaction network perspective.  The label of the reaction indicates also the chosen reaction rate. The nonnegative parameters $\lambda,\mu_s, \mu_i, \mu_r$ are the birth rate of $S$ and the death rates of $S,I,R$ respectively:
\begin{equation*}
... \quad \underset{\lambda}{\longrightarrow}\quad S, \quad\quad\quad
S \quad \underset{\mu_s s}{\longrightarrow} \quad...\;,\quad\quad\quad
I \quad \underset{\mu_i i}{\longrightarrow} \quad...\;,\quad\quad\quad
R \quad \underset{\mu_r r}{\longrightarrow} \quad...\; .
\end{equation*}
Standard linear reactions are assumed to model vaccinations, the transition from being recovered to being susceptible, and and the unproblematic transition from being infected to being recovered:
\begin{equation*}
S \quad \underset{\gamma_s s}{\longrightarrow} \quad R, \quad\quad\quad
R \quad \underset{\gamma_r r}{\longrightarrow} \quad S,\quad\quad\quad I \quad \underset{\gamma_i i}{\longrightarrow} \quad R.\end{equation*}
The only nonlinearities appear in the infection rate and potentially in the \emph{treatment rate}
\begin{equation*}S+I \quad \underset{\beta s i F(i)}{\longrightarrow} \quad 2I,
\quad\quad\quad I \quad \underset{T(i)}{\longrightarrow} \quad R.\end{equation*}

\paragraph{Historical account}
The SIRS model \eqref{Vys} generalizes several papers 
\cite{Capasso,LiuLevin,WR,Wang,ZhouFan,Vyska,LuRuan19,xu2021complex,Roostaei}. This epidemic model is firstly inspired by the pioneering work of Capasso and Serio \cite{Capasso}, which proposed to capture in $F(i)$ self-regulating adaptations of the population to the epidemics due to the spread of information on its existence. At the origin $F(i)$ must be $1$ to fit the Kermack-McKendrick SIR model, and the authors took  $F(i) =\frac{1}{1 + b i}$, where $b$ measures the  \emph{inhibition effect}. The ensuing infection rate $f(s,i):=\beta s i F(i)$ then follows Michaelis--Menten in $i$.

A second strain of literature was initiated by Wang and Ruan \cite{WR,Wang} (see also
 \cite{Vyska}), who included a discrete treatment term $T(i)=\eta \operatorname{min}[i, w]$   which takes into account that some people recover only after having been treated in hospital. The minimum captures
 a possible breakdown of the medical system due to the epidemics. Subsequently, Zhang and Liu \cite{zhang2008backward} proposed to study  smooth treatment terms $T(i)=\eta \frac{i}{1+ \nu i}$.
In both versions
  $\eta=T'(0)$ is  the treatment rate when $i=0$,
and  $w=\frac{\eta}{\nu}$ is an upper limit for the capacity of treatment. Roostaei et al. \cite{Roostaei} further simplified the model by assuming $\eta=1$ and showed the occurrence of Hopf bifurcations. Putting together the nonlinear term of Capasso with that of Wang and Ruan gave rise to a huge literature focused on understanding the possible bifurcations which may arise. Further explorations concerned different nonmonotone nonlinearities for $F(i)$  \cite{Ruan07, ZhouXiao, LuRuan19, Zhang23}.


We recall again a fundamental result of ME: the \emph{next generation matrix theorem} \cite{Van08, Van}. This result states that the DFE of an epidemic model
is locally stable if the basic reproduction number $R_0< 1$  and
 unstable if $R_0> 1$. Let $\bar{s}_{DFE}$ indicate the value of $s$ at the DFE. It may be shown  that the $R_0$ for the  model \eqref{Vys}   is:
 \begin{equation}\label{R0}
 R_0:= \bar{s}_{DFE} \mathcal{R}\quad \text{where}\quad\mathcal{R}:=\frac{\beta F(0)}{\mu_i +\gamma_i   +T'(0)}.
 \end{equation}
  
 Particular cases of this formula have appeared in the literature \cite{ZhouFan, Vyska, Gupta}; however,
$\bar{s}_{DFE}$  is absent from \eqref{R0} in these papers, as they study a particular case with  $\gamma_s=0$, which implies $\bar{s}_{DFE}=1$. Note that \eqref{R0} suggests that models in which $F(0)=0$ should perhaps not even be accepted as `epidemiological', since $R_0$ does not have the usual epidemiological interpretation. One of such examples \cite{LuRuan19} is $ F(i)= \frac{i}{1+ b_1 i + b_2 i^2}$ and thus $F(0)=0$. Interestingly, this model shows a variety of dynamical behaviors as consequence of \emph{cusp} and \emph{Takens-Bogdanov bifurcations}. In later work\cite{Zhang23}, Zhang and Li showed that $F(i) =\frac{1+a i}{1  + b i^2}$ leads to similar features, with the significant improvement that $R_0$ may be bigger than $1$.

\subsection{General results}

For convenience, we define
\begin{equation}
\begin{cases}
f(s,i):=\beta s i F(i) \quad\text{ with $F(i)> 0, \;F(0)=1$};\\
h(i):=\gamma_i i+T(i) \quad\text{  with $\gamma_i>0$, $T(i)\ge0$, $T'(i)\ge0$.}\
\end{cases}
\end{equation}

It is easy to see that an endemic equilibrium $(\bar{s},\bar{i},\bar{r})>0$ characterizes the flux constraints:
\begin{equation}\label{eq:endeqcap}
\begin{cases}
\lambda=\mu_s\bar{s}+\mu_i\bar{i}+\mu_r \bar{r};\\
f(\bar{i})=h(\bar{i})+\mu_i \bar{i};\\
\gamma_r\bar{r}+\mu_r \bar{r}=\gamma_s\bar{s}+h(\bar{i}).
\end{cases}
\end{equation}
The constraints \eqref{eq:endeqcap} together with $h(i)>0$ imply that the mere existence of an endemic equilibrium necessarily requires 
\begin{equation}\label{eq:capassoeqnecessary}\gamma_r+\mu_r>0.
\end{equation}

The symbolic Jacobian of the system \eqref{Vys} then reads:
\begin{equation}\label{eq:jac}
G=\begin{pmatrix}
    -f'_s - \gamma_s-\mu_s & -f'_i & \gamma_r\\
    f'_s & f'_i - h' -\mu_i & 0\\
    \gamma_s & h' & -\gamma_r -\mu_r
\end{pmatrix}.
\end{equation}

\begin{obs}\label{obs:gersh}
If $f'_i\le 0$ 
 then the symbolic Jacobian matrix $G$ is weakly-diagonally dominant, i.e. $|G_{ii}|\ge \sum_j |G_{ji}|$, and thus $G$ never possesses eigenvalues with positive real-part. This excludes any bifurcation at endemic equilibria.
\end{obs}
 In particular, a bifurcation requires that $f'_i>0$ at an endemic equilibrium. We assume this throughout this section. Our main result of this section strengthens Observation \ref{obs:gersh}.

\begin{theorem}\label{thm:necessarybif}
For any bifurcation to occur at an  endemic equilibrium  $(\bar{s},\bar{i},\bar{r})$ of system \eqref{Vys} it is necessary that
\begin{equation}
G_{ii}(\bar{s},\bar{i})>0\end{equation}
or - equivalently -
\begin{equation}\label{eq:Rcapasso}
f'_i(\bar{s},\bar{i})>h'(\bar{i})+\mu_i.
\end{equation}
In particular, either $f$ or $h$ must be nonlinear in $i$.
\end{theorem}

\begin{remark} \label{rq:Rsi}
Theorem \ref{thm:necessaryHopfthm} concerns endemic equilibria. However, \eqref{eq:Rcapasso} generalizes \eqref{R0}:
$$f'_i(\bar{s},\bar{i})=\beta \bar{s}F(\bar{i})+\beta\bar{s}\bar{i}F'(\bar{i}),$$
and thus the condition for DFE equilibria $(\bar{s}_{DFE}, 0,\bar{r})$ reads precisely as in \eqref{R0}. 
\end{remark}

As a special but important case consider
\begin{equation}\label{eq:rationalform}
f(s,i)=\beta s i\frac{1}{1+b i^n}, \quad n\ge 1,
\end{equation}
which includes in particular Michaelis-Menten for $n=1$, and the case of (nonmonotone!) \emph{generalized Monod-Haldane} function (see for example \cite{sokol81} when $n=2$). We have the following Corollary.

\begin{cor}\label{cor:necessaryhopfMM}
If $f(s,i)$ belongs to the Monod-Haldane class \eqref{eq:rationalform}, then for an equilibrium bifurcation to occur, the nonlinearity of $T(i)$ is necessary.
\end{cor}

\proof
Firstly note that
\begin{equation}
f'_i(s,i)=\frac{\beta s(1-b(n-1) i^n)}{(1+b i^n)^2} \le \frac{\beta s}{(1+b i^n)^2}=\frac{\beta s i}{(1+b i^n)^2 i}=\frac{f(s,i)}{(1+b i^n)i}\le\frac{f(s,i)}{i}.
\end{equation}
Assume  that $h(i)$ is linear, i.e. $T(i)=0$. Fix an equilibrium $(\bar{s},\bar{i})$. We have
\begin{equation}
G_{ii}(\bar{s},\bar{i})=f'_i(\bar{s},\bar{i})-\gamma_i-\mu_i\le\frac{f(\bar{s},\bar{i})}{\bar{i}}-\gamma_i -\mu_i =\frac{f(\bar{s},\bar{i})-\gamma_i \bar{i}  -\mu_i \bar{i}}{\bar{i}}=0
\end{equation}
and Theorem \ref{thm:necessarybif}
excludes equilibria bifurcations.
\endproof

Instead of proving directly Theorem \ref{thm:necessarybif}, we split the problem into \emph{zero-eigenvalue bifurcations} and \emph{Hopf} (or purely-imaginary) bifurcations, and prove the results in Lemmas \ref{lem:necessarysn} and \ref{thm:necessaryHopfthm} below. Both lemmas jointly constitute Theorem \ref{thm:necessarybif}. However. As proofs and applications are quite different, we proceed in parallel.

\paragraph{Zero-eigenvalue bifurcations.}

First, it is worth noting that
there are models \eqref{Vys} that do not admit zero-eigenvalue bifurcations for any choice of infection rate. To see this, we compute the characteristic polynomial of the symbolic Jacobian $G$:
$$g(\lambda)=\operatorname{det}(G-\lambda \operatorname{Id})=c_0\lambda^3+c_1\lambda^2+c_2\lambda+c_3,$$
\begin{equation}\label{eq:cp_G}
\text{with }
\begin{cases}
\begin{split}
c_0&=-1\\
c_1&=f'_i - f'_s - \gamma_r - \gamma_s - h' - \mu_s-\mu_i - \mu_r\\
c_2&=f'_i(\gamma_r + \gamma_s + \mu_r + \mu_s) - f'_s\gamma_r - f'_sh' - \gamma_rh' - \gamma_sh' - f'_s\mu_i - f'_s\mu_r \\
&- \gamma_r\mu_i - \gamma_s\mu_i - \gamma_r\mu_s - \gamma_s\mu_r - h'\mu_r - h'\mu_s - \mu_i\mu_r - \mu_i \mu_s - \mu_r \mu_s\\
c_3&= f'_i (\mu_s \gamma_r+\mu_s\mu_r+\gamma_s\mu_r)-\gamma_sh'\mu_r-\gamma_s\mu_i\mu_r-\mu_s  h' \gamma_r -f'_s h' \mu_r -  \mu_i \mu_s \gamma_r\\ &- f'_s \mu_i \mu_r-f'_s \mu_i\gamma_r -\mu_s h' \mu_r-\mu_s \mu_i \mu_r,
\end{split}
\end{cases}
\end{equation}
From a simple inspection of \eqref{eq:cp_G} we derive the following straightforward lemma.
\begin{lemma}\label{lem:zeroeigcap}
The following statements hold:
\begin{enumerate}
    \item For the fully-open model $\mu_s,\mu_i,\mu_r>0$ and the fully-closed model $\mu_s=\mu_i=\mu_r=0$
there is always a choice of $f'_i>0$ such that the symbolic Jacobian $G$ has a change of stability via zero-eigenvalue.
\item Consider the half-open models with $\mu_i>0$, $\mu_s=0$. If  $\gamma_s\mu_r=0$ then there is never a choice of parameters such that the symbolic Jacobian $G$ has a change of stability via zero-eigenvalue.
\end{enumerate}
\end{lemma}
\proof
Firstly recall that the necessary condition \eqref{eq:capassoeqnecessary} for the existence of endemic equilibria. The proof follows from \eqref{eq:cp_G}. If $\mu_s,\mu_i,\mu_r>0$, $c_3$ (the determinant of $G$) is a multilinear homogenous polynomial with monomials of opposite sign. In particular, $\mu_sf'_i\mu_r$ has a positive coefficient and $\mu_s\mu_i\mu_r$ has a negative coefficient. Since the monomials where $f'_i$ appears identify the monomials with a positive coefficient, we can always solve for $f'_i>0$.\\ 
Analogously, if $\mu_s=\mu_i=\mu_r=0$ then $c_3\equiv 0$. There is indeed a conservation law (the total population) and a thus a trivial eigenvalue zero. However, as before, $c_2$ presents terms with opposite signs where the positive terms are characterized by $f'_i$. Indeed, the bare existence of an endemic equilibrium with $\mu_r=0$ requires $\gamma_r>0$, via \eqref{eq:capassoeqnecessary}. Then $c_2$ has always at least two nonzero summands, $f'_i\gamma_r$ and  $f'_s\gamma_r$, with opposite sign, and we can solve $c_2=0$ for $f'_i>0$ and find parameter values where the Symbolic Jacobian has a non-trivial zero-eigenvalue.\\
On the other hand, if $\mu_s=\gamma_s\mu_r=0,$  we get that $c_3\le0$. Moreover,  $\mu_i>0$ and \eqref{eq:capassoeqnecessary} imply  $-f'_s\mu_i(\gamma_r+\mu_r)<0$ and thus $c_3<0$.
\endproof

The half-open models in point 2 of Lemma \ref{lem:zeroeigcap} can be thought to model deadly infectious diseases where the death rate of susceptible population is considered neglectable.

Besides general symbolic considerations, in practice given epidemiological models within \eqref{Vys} assume a rather specific choice in the reaction rates, and it is thus important to have an idea of which choice potentially leads to bifurcations. We have the following result.

\begin{lemma}\label{lem:necessarysn}
For a zero-eigenvalue bifurcation to occur at an endemic equilibrium point $(\bar{s},\bar{i},\bar{r})$ of system $\eqref{Vys}$ it is necessary that
\begin{equation}
G_{ii}(\bar{s},\bar{i})>0\end{equation}
or - equivalently -
\begin{equation}
f'_i(\bar{s},\bar{i})>h'(\bar{i})+\mu_i.
\end{equation}
In particular, either $f$ or $h$ must be nonlinear in $i$.
\end{lemma}
\proof
The proof again follows simply by inspection of the characteristic polynomial \eqref{eq:cp_G}. We rewrite the coefficients $c_2, c_3$ as follows
\begin{equation*}
\begin{cases}
c_2=(f'_1-h'-\mu_i)(\gamma_r+\gamma_s+\mu_r+\mu_s)-f'_s(h'+\mu_i+\gamma_r+\mu_r)-\mu_s(h'+\gamma_r+\mu_r)\\
c_3=(f'_i-h'-\mu_i) (\mu_s \gamma_r+\mu_s\mu_r+\gamma_s\mu_r) -f'_s (h' \mu_r +   \mu_i \mu_r+\mu_i\gamma_r)
\end{cases},
\end{equation*}
from which it is clear that if $G_{ii}=f'_1-h'_1-\mu_i<0$ then the characteristic polynomial has all strictly negative coefficients and thus no zero root, i.e. no zero-eigenvalues are possible. Alternatively, note that Descartes' Rule of sign excludes real positive roots.

Finally, indirectly assume that $f$ and $g$ are both linear in $i$, i.e. $F(i)\equiv 1$, $T(i)\equiv 0$ and $f=\beta s i $, $h(i)=\gamma_i i$. At an endemic equilibrium value, we get
\begin{equation}\beta \bar{s} \bar{i} - \gamma_i \bar{i}  - \mu_i \bar{i}=0,\end{equation}
which implies \begin{equation}
 \beta \bar{s}=\gamma_i+\mu_i,
\end{equation}
contradicting $f'_i=\beta\bar{s}>h'+\mu_i=\gamma_i+\mu_i$;  hence either  $f$ or $h$ must be nonlinear in $i$.
\endproof

\paragraph{Hopf bifurcations.}

Firstly, from a symbolic perspective, Theorem \ref{lem:purelyimaginary} applies and guarantees the following corollary:
\begin{cor}\label{cor:hopfhalfclosed}
There is a choice of positive $f'_s,  f'_i, h', \gamma_s, \mu_s, \mu_i, \gamma_r, \mu_r>0$ for which the Symbolic Jacobian $G$, \eqref{eq:jac}, has purely imaginary eigenvalues.
\end{cor}

However, the more realistic reaction rates of system \eqref{Vys} involve many linear terms in mass action form, and thus it is not guaranteed that we can always find Hopf bifurcation, irrespective of the functional nonlinear form of $f(s,i)$ and $T(i)$. We have the following result.

\begin{lemma}\label{thm:necessaryHopfthm}
For a Hopf bifurcation at an endemic equilibrium point $(\bar{s},\bar{i},\bar{r})>0$, it is necessary that
\begin{equation}
G_{ii}(\bar{s},\bar{i})>0,
\end{equation}
or - equivalently -
\begin{equation}
f'_i(\bar{s},\bar{i})>h'(\bar{i})+\mu_i.
\end{equation}
In particular, either $f$ or $h$ must be nonlinear in $i$.
\end{lemma}

\begin{remark} Lemma \ref{thm:necessaryHopfthm} strengthens the well-known result that the purely mass-action SIRS does not have Hopf points.
\end{remark}

\proof
Firstly, we show that the matrix $A$ defined as
\begin{equation}\label{A}
A=\begin{pmatrix}
    -f'_s - \gamma_s-\mu_s & -f'_i & \gamma_r\\
    f'_s & -\alpha & 0\\
    \gamma_s & h' & -\gamma_r -\mu_r
\end{pmatrix},
\end{equation}
does not admit purely imaginary eigenvalues for any choice of $$f'_s, f'_i, h'>0 \quad \text{and} \quad \gamma_s, \mu_s, \gamma_r, \mu_r, \alpha \ge0.$$

Indeed: consider the characteristic polynomial of such $A$, evaluated at purely imaginary $i\omega$. Here $i$ is the imaginary unit, of course.

\begin{equation*}
\begin{split}
P(A):&=\operatorname{det}(A-\omega i \operatorname{Id})\\ &= \omega^3i+ (f'_s + \gamma_s+ \mu_s + \alpha + \gamma_r  + \mu_r)\omega^2\\ &
+(- \alpha f'_s - \alpha \gamma_r - \alpha \gamma_s - f'_if'_s - f'_s\gamma_r - \alpha\mu_r - \alpha\mu_s - f'_s\mu_r - \gamma_r\mu_s - \gamma_s\mu_r - \mu_r\mu_s) \omega i\\ &-
 f'_s\alpha\gamma_r - f'_sf'_i\gamma_r - f'_s\alpha  \mu_r - \mu_s\alpha \gamma_r  - \gamma_s\alpha \mu_r + f'_s h' \gamma_r - f'_sf'_i\mu_r - \mu_s\alpha\mu_r
\end{split}
\end{equation*}
Split $P(A)$ in imaginary and real part:
\begin{equation*}
\begin{split}
\Re(P(A))&=(f'_s + \gamma_s+ \mu_s + \alpha + \gamma_r  + \mu_r)\omega^2--
 f'_s\alpha\gamma_r - f'_sf'_i\gamma_r - f'_s\alpha  \mu_r - \mu_s\alpha \gamma_r  \\&- \gamma_s\alpha \mu_r + f'_s h' \gamma_r - f'_sf'_i\mu_r - \mu_s\alpha\mu_r\\
\Im(P(A))&=\omega(\omega^2 - \alpha f'_s - \alpha \gamma_r - \alpha \gamma_s - f'_if'_s - f'_s\gamma_r - \alpha\mu_r - \alpha\mu_s - f'_s\mu_r - \gamma_r\mu_s - \gamma_s\mu_r - \mu_r\mu_s)
\end{split}
\end{equation*}
We can divide $\Im(P(A))$ by $\omega$, eliminating the trivial real root $\omega=0$:
\begin{equation}\tilde{\Im}(P(G))=\Im(P(G))/\omega
\end{equation}
Define $\Omega=\omega^2$. Purely imaginary eigenvalues necessarily require
\begin{equation}\operatorname{resultant}(\Re(P(G)),\tilde{\Im}(P(G)),\Omega)=0.
\end{equation}

However, a direct computation shows \begin{equation}
\begin{split}
    \operatorname{resultant}(\Re(P(G),\Im(P(G),\Omega)&=- f'_s\alpha^2 - \alpha^2\gamma_r - \alpha^2 \gamma_s - \alpha^2 \mu_r - \alpha^2 \mu_s - \alpha  {f'_s}^2\\& - 2 \alpha f'_s \gamma_r - 2 \alpha f'_s \gamma_s - 2 \alpha f'_s \mu_r - 2 \alpha f'_s \mu_s - f'_i \alpha f'_s\\ &- \alpha \gamma_r^2 - 2 \alpha \gamma_r \gamma_s - 2 \alpha \gamma_r \mu_r - 2 \alpha \gamma_r \mu_s - \alpha \gamma_s^2 \\ &- 2 \alpha \gamma_s \mu_r - 2 \alpha \gamma_s \mu_s - \alpha \mu_r^2 - 2 \alpha \mu_r \mu_s - \alpha \mu_s^2\\ &-  {f'_s}^2 \gamma_r -  {f'_s}^2 \mu_r - f'_i  {f'_s}^2 - f'_s \gamma_r^2 - f'_s \gamma_r \gamma_s \\&- 2 f'_s \gamma_r \mu_r - 2 f'_s \gamma_r \mu_s - h' f'_s \gamma_r - 2 f'_s \gamma_s \mu_r \\
    &- f'_i f'_s \gamma_s - f'_s \mu_r^2 - 2 f'_s \mu_r \mu_s - f'_i f'_s \mu_s\\ &- \gamma_r^2 \mu_s - \gamma_r \gamma_s \mu_r - \gamma_r \gamma_s \mu_s - 2 \gamma_r \mu_r \mu_s\\ &- \gamma_r \mu_s^2 - \gamma_s^2 \mu_r - \gamma_s \mu_r^2 - 2 \gamma_s \mu_r \mu_s - \mu_r^2 \mu_s - \mu_r \mu_s^2\\
 &\neq 0
\end{split}
\end{equation}
for any positive choice of symbols. The presence of monomials ${f'_s} \alpha$, ${f'_s}^2\gamma_r$, ${f'_s}^2\mu_r$, $f'_s\gamma_sf'_i$, $f'_s\mu_sf'_i$  guarantees the resultant being nonzero for any choice of  $\gamma_s, \mu_s, \gamma_r, \mu_r, \alpha$ where at least one of the symbols is nonzero. To conclude, consider the case $$\gamma_s, \mu_s, \gamma_r, \mu_r, \alpha=0.$$
The matrix $A$ in \eqref{A} becomes
\begin{equation}
A(0)=\begin{pmatrix}
    -f'_s & -f'_i & 0\\
    f'_s & 0 & 0\\
    0 & h' & 0
\end{pmatrix},
\end{equation}
which has one eigenvalue $0$ and $\operatorname{tr}(A)=-f'_s\neq 0$: in particular, $A(0)$ does not admit purely imaginary eigenvalues.

To prove Lemma \ref{thm:necessaryHopfthm}, it is enough to notice that the set of matrices $A$, \eqref{A}, strictly contains any choice of symbolic Jacobian $G$ with nonpositive $i$-diagonal entry
\begin{equation}
    f'_i-h_i'-\mu_i\le 0.
\end{equation}
The nonlinearity of either $f$ or $g$ follows from $G_{ii}\neq 0$, just as in Proof of Lemma \ref{lem:necessarysn}.
\endproof

\subsection{Explicit examples.}\label{sec:explicitexC}

\paragraph{Absence of bifurcations.} Observation \ref{obs:gersh} implies that no bifurcation is possible at an endemic equilibrium with a decreasing infection rate. On the other hand, the original model of Capasso and Serio \cite{Capasso} reads:
\begin{equation}\label{eq:originalcapasso}
\begin{cases}
\dot{s}= -\gamma_s s - f(s,i),\\
i'=f(s,i)-\gamma_ii,\\
\dot{r}=\gamma_s s+ \gamma_ii,
\end{cases}
\end{equation}
where $\gamma_s \ge 0$ and $f(s,i)=\beta s i F(i)$ is
in Michaelis-Menten form in $i$: $$f(s,i)=\beta s i \frac{1}{1+b i},$$
with  inhibition rate $F(i)= \frac{1}{1+ b i}$.
The original model \eqref{eq:originalcapasso} does not possess any endemic equilibrium. Therefore, we may consider variants of \eqref{eq:originalcapasso},
where we either add demography ($\lambda, \mu_s,\mu_i,\mu_r > 0$),
 or a transition rate from the recovered state to a newly susceptible state
 ($\gamma_r\neq 0$). The resulting open model where we add both is just
\eqref{Vys} with linear $h(i)=\gamma_i i$, and the closed version
 is \eqref{Vys} with linear $h(i)=\gamma_i i$ and   $\lambda,\mu_s,\mu_i\mu_r=0$. Corollary \ref{cor:necessaryhopfMM}
  implies that none of these extensions admit any equilibria bifurcation: the endemic point is never unstable: either a
  different non-linearity in $f$ or a nonlinear treatment term $T(i)$ is needed.

  Similarly, no bifurcation occurs for the model of Xuan and Ruan \cite{Ruan07}, as they considered linear zero treatment $T(i)=0$, and consequently a linear $h(i)=\gamma_i i$. The infection rate is in Monod-Haldane form $$f(s,i)=\beta s \frac{i}{1+b i^2}.$$

Again, Corollary \ref{cor:necessaryhopfMM} excludes bifurcations.

\paragraph{Existence of zero-eigenvalue bifurcation.}
We show the existence of zero-eigenvalue bifurcations with a minimal presence of nonlinearity in the system.
\begin{lemma}\label{lem:zeroeigRW}
For the Capasso-Ruan-Wang model \eqref{Vys} with constant inhibition function $$F(i)=1,$$ 
and treatment rate in Michaelis-Menten form:
\begin{equation}
    T(i)= a\frac{i}{1+ b i},
\end{equation}
there exist endemic equilibria $(\bar{s},\bar{i},\bar{r})>0$ with a zero-eigenvalue Jacobian.
\end{lemma}

The choice of constant inhibition function reduces the infection rate $f(s,i)=\beta s i$ in mass action form. Treatment terms in Michaelis-Menten form have been studied in the literature \cite{zhang2008backward,Roostaei}, where however the system reduced eventually to a 2-dimensional  problems. Since we allowed also for loss of immunity via $\gamma_r$, our system may remain three-dimensional.

\proof We follow Section \ref{sec:sna}: We firstly observe that the equilibria constraints \eqref{eq:endeqcap} lead to a natural parametrization of the fluxes at the endemic equilibrium $(\bar{s},\bar{i},\bar{r})>0$. Let $d_1,d_2,d_2,\theta, g>0$ be positive parameters. Then the fluxes at an endemic equilibrium of \eqref{Vys} are parametrized as:
\begin{equation}
\begin{pmatrix}
    \lambda\\
    \beta \bar{s} \bar{i}\\
    h(\bar{i})\\
    \gamma_s \bar{s}\\
    \gamma_r \bar{r}\\
    \mu_s \bar{s}\\
    \mu_i \bar{i}\\
    \mu_r \bar{r}
\end{pmatrix}=
\begin{pmatrix}
    d_1+d_2+d_3\\
    \theta+d_2\\
    \theta\\
    g\\
   g+\theta-d_3\\
    d_1\\
    d_2\\
    d_3
\end{pmatrix}.
\end{equation}
Note that $h(i)=\gamma_i i + T(i)$ satisfies:
$$h'(i)=\gamma_i+\frac{a}{(1+bi)^2}=\bigg(\gamma_i i + \frac{T(i)}{(1+bi)}\bigg)\frac{1}{i},$$
so that for proper choices of $\gamma_i, b$, the map $h'(\bar{i})(b,\gamma_i)$ is surjective onto $(0,\theta/\bar{i})$. Following our Theorem \ref{thm:paramquasi}, we then define
$$q=q(b,\gamma_i)\in [1,+\infty) \quad \text{and, for simplicity,} \quad \vartheta:=\theta/q\in (0,\theta]$$
As a consequence the Jacobian matrix
\begin{equation}
Jac=\begin{pmatrix}
-\beta i -\gamma_s -\mu_s &-\beta s & \gamma_r\\
\beta i & \beta_s-h'-\mu_i & 0\\
\gamma_s & h' & -\gamma_r-\mu_r
\end{pmatrix}
\end{equation}
at $(\bar{s},\bar{i},\bar{r})$ is parametrized as follows:
\begin{equation}
\begin{split}
Jac=&Jac \;\operatorname{Id}=Jac \;\operatorname{diag}(\bar{s},\bar{i},\bar{r})\;\operatorname{diag}(1/\bar{s},1/\bar{i},1/\bar{r})\\
=&\begin{pmatrix}
-\theta - d_1-d_2 -g & -\theta-d_2 & g+\theta-d_3\\
\theta+d_2 & \theta-\vartheta & 0\\
g & \vartheta & -g-\theta
\end{pmatrix}\operatorname{diag}(1/\bar{s},1/\bar{i},1/\bar{r})\\
&=B(v)\operatorname{diag}(1/\bar{s},1/\bar{i},1/\bar{r})
\end{split}
\end{equation}
Note again that
$$\operatorname{rank}B(v)\operatorname{diag}(1/\bar{s},1/\bar{i},1/\bar{r})=\operatorname{rank}B(v),$$
and thus $B(v)$ fully identify zero-eigenvalue bifurcations. The determinant of $B(v)$ reads:
$$\operatorname{det}Jac=(d_1 - d_2)\theta^2 + (- d_2^2 - gd_2 - \vartheta d_1 - \vartheta d_3 + d_1 g + d_3 g)\theta - gd_2^2 - \vartheta d_3d_2 - \vartheta d_1 g - \vartheta d_3 g.$$
It is straightforward to find several instances for a zero-eigenvalue at endemic equilibria, solving for $\theta$. For instance, if $d_1>d_2$ as a consequence of Descartes' rule of sign, or in the case $d_1=d_2=d_3:=d$ if
$$dg - 2\vartheta d - d^2>0,$$
which definitely holds for $g$ sufficiently large.
\endproof

\begin{remark} Numerical simulations indicate indeed that the zero eigenvalue bifurcations above are of saddle-node type. \end{remark} 

As a final consideration, for the closed system without demographics, i.e. $d_1=d_2=d_3=0$, we get the following characteristic polynomial of $Jac_{c}$:
$$P(Jac_c,\lambda)=- \lambda^3 + (- a - 2g - h)\lambda^2 + (g\theta - 2\vartheta \theta - 2\vartheta)\lambda,$$
which admits nontrivial zero eigenvalues only if and only if (!) $g$ is sufficiently large. In particular, if $\gamma_s=g=0$, i.e., no vaccinations, then the closed system does not admit zero-eigenvalue bifurcation even with nonlinear treatment.

\paragraph{Existence of Hopf bifurcation.}
Theorem \ref{thm:necessaryHopfthm} and Corollary \ref{cor:necessaryhopfMM} imply that $T(i)$ must be nonlinear for Hopf bifurcations in systems where $f(s,i)$ belongs to the Monod-Haldane class. As we did for zero-eigenvalue bifurcations, we show here that a nonlinear $T(i)$ is sufficient to have Hopf bifurcations, even if all the other rates are in mass-action form. More specifically, we hunt Hopf bifurcations again for the simplest case: $f(s,i)=\beta s i$ and a Michaelis-Menten treatment:
\begin{equation}
    T(i)= a\frac{i}{1+ b i}.
\end{equation}



The convenience of the symbolic approach discussed in Section \ref{sec:symbolichunt} is that we do not start by computing explicitly the equilibria, and then looking for instances of the Jacobian at such equilibria with purely imaginary eigenvalues: such approach renders computations nontrivial. Instead, we first spot in the symbolic Jacobian \eqref{eq:jac} the derivatives that might trigger bifurcations.
In this heuristic sense, Theorem \ref{lem:purelyimaginary} suggests the natural choice of bifurcation parameter  $\xi:=f'_i(\bar{s},\bar{i})=\beta \bar{s}$, where the bifurcating equilibrium $(\bar{s},\bar{i}, \bar{r})$ is yet-to-be defined.

Instead of venturing into confusing symbolic computations, we focus on one single parameter-choice, for sake of explanation. We arbitrarily choose simple values for
the other derivatives besides $\xi$, which will also be used later in bifurcation plots:
\begin{equation}\label{eq:partialderex}
\begin{cases}
f'_s=\gamma_s=\gamma_r=\mu_r=\gamma_i=T'_i=1;\\
\mu_s=0.1;\\
h'=\mu_i=2.\\
\end{cases}
\end{equation}

With the derivative values above,   the 1-parameter  symbolic Jacobian is:
\begin{equation}\label{eq:Gxi}
G(\xi )=\begin{pmatrix}
    -1 -1 - 0.1 & -\xi & 1\\
    1 & \xi  -2-2 & 0\\
    1 & 2 & -1-1
\end{pmatrix}.
\end{equation}

  The characteristic polynomial corresponding to \eqref{eq:Gxi} evaluated at a purely imaginary eigenvalue is given by
$$ P(G(\xi))=Det[G(\xi)-  \; \omega i \;  \text{Id} ],$$
with real and imaginary parts given, respectively, by
\begin{equation}
\begin{cases}
Re[P(G)]=-\xi  \omega ^2+\frac{6 \xi }{5}+\frac{81 \omega ^2}{10}-\frac{54}{5},\\
Im[P(G)]=\omega ^3+\frac{31 \xi  \omega }{10}-\frac{98 \omega }{5}
\end{cases}
\end{equation}
After dividing $Im[P(G)]$  by the trivial root $\omega$, we obtain $$\Tilde{Im[P(G)]}:= Im[P(G)]/{\omega}=\frac{31 \xi }{10}+\omega ^2-\frac{98}{5},$$
and putting $\Omega=\omega^2$,  we obtain the resultant:
$$\text{Resultant}[Re[P(G)],\Tilde{Im[P(G)]},\Omega]=\frac{31 \xi ^2}{10}-\frac{4351 \xi }{100}+\frac{3699}{25}.$$
The last equation equals zero when
$$\xi=\xi^*_{-,+}:=\frac{1}{620} \left(4351\pm \sqrt{584161}\right)=(5.78499,8.25049).$$

As a sanity check,
  the eigenvalues of $G(\xi^*_-)$ are indeed $(-2.31501,\pm 1.29094 \; Im)$: a purely imaginary crossing.

  We turn now to the parametrization of the nonlinear treatment function $T(i)$, so that the matrix $G(\xi)$ is the actual Jacobian of the system, for a $\xi$-parametrized family of equilibria. Before that, we must specify all the parameters, and one simplest choice for specifying the first derivative in \eqref{eq:partialderex} is  $\bar{i}=1, \beta=1$  and hence $\xi =\bar{s}$.

  With this choice, the nonlinear $h(i)=\gamma_i i + T(i)$ must satisfy at the endemic  equilibrium
\begin{equation}
\begin{cases}
h(1)=1+T(1)=\xi - 2\\
h'(1)=1 + T'(1)=2,
\end{cases}
\end{equation}
that is
\begin{equation}
\begin{cases}
T(1)=\xi - 3\\
T'(1)=1,
\end{cases}
\end{equation}
With $T(i)$ in Michaelis-Menten form
$$T(i)=a\frac{i}{1+bi},$$
we can solve for $a,b>0$ using the equality
\begin{equation}
T'(i)=\frac{a}{(1+bi)^2}=\frac{T(i)}{(1+bi)i}.
\end{equation}
At $i=1$ we get
\begin{equation}
\begin{cases}
b=\frac{T(1)}{T'(1)}-1=\xi-4\\
a=T(1)(1+b)=(\xi-3)(1+\xi-4)=(\xi-3)^2\\
\end{cases}
\end{equation}
This leads to
$$T(i)=(\xi-3)^2\frac{i}{1+(\xi-4)i}.$$
Looking at the $r$ equilibrium equation, we get
$$0=\xi+1+\xi-3-\gamma_r r-\mu_r r.$$
Together with the chosen Jacobian \eqref{eq:Gxi}, this implies
\begin{equation}
\bar{r}=\xi -1.
\end{equation}
Finally, solving for $\lambda$ the equilibrium equation for $s$ we get
$$\lambda=\xi+\xi-\xi+1+0.1\xi=(1.1\xi+1).$$
The chosen reaction rates give rise to the following 1-parameter version $C(\xi )$ of the model \eqref{Vys}:
\begin{equation}\label{sirxi}
\begin{cases}
\dot{s}=(1.1\xi +1) -s-si+ r - 0.1 s\\
i'=si-i-(\xi -3)^2\frac{i}{1+(\xi -4)i}-2 i\\
\dot{r}=s+i+(\xi -3)^2\frac{i}{1+(\xi -4)i}- r-r
\end{cases},
\end{equation}
with 1-parameter family of equilibria $(\bar{s}(\xi ),\bar{i}(\xi ),\bar{r}(\xi ))=(\xi ,1,\xi -1)$. By construction, the Jacobian \eqref{eq:jac} evaluated at an equilibrium $(\bar{s}(\xi ),\bar{i}(\xi ),\bar{r}(\xi ))$ of $C(\xi )$ reads as \eqref{eq:Gxi}, and thus the system undergoes a Hopf bifurcation for values $\xi=\xi^*_-=\frac{1}{620} (4351- \sqrt{584161})$.

We confirmed numerically that for $\xi \geq \xi_H$  there exists a stable limit cycle around the equilibrium $(\xi,1,\xi-1)$,  see Figures \ref{fig:1a}, \ref{fig:1b}, \ref{fig:1c} below.

    \begin{figure}[hbt!]
    \centering
\includegraphics[width=0.5\textwidth]{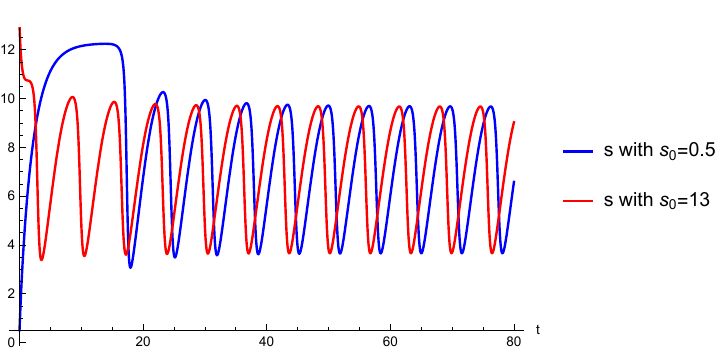}
        \caption{\footnotesize $(s)$-time plot at $\xi_H=\frac{1}{620} (4351- \sqrt{584161})$ (when the eigenvalues of $G(\xi)$ are $(-2.31501,\pm 1.29094 \; Im)$) with two different initial values $(0.5,13)$, suggest the existence of a unique limit stable cycle.}
        \label{fig:1a}
    \end{figure}%
    ~
     \begin{figure}[hbt!]
     \centering
\includegraphics[width=0.5\textwidth]{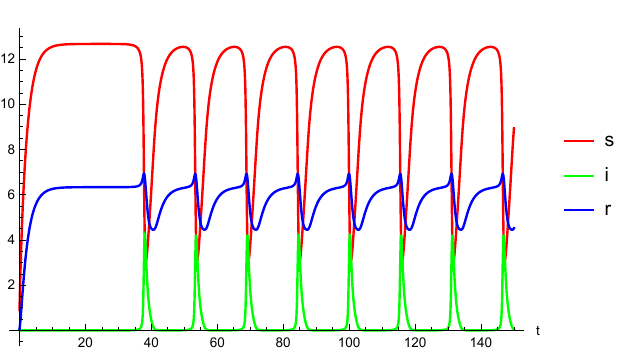}
        \caption{\footnotesize $(s,i,r)$-time plot indicate the existence of periodic solutions when $\xi=6$. In this case the eigenvalues of $G(\xi)$ at the (unstable fixed point) are $(-2.33284,0.116421 \pm 1.23678\; Im)$}
        \label{fig:1b}
    \end{figure}
     ~
     \begin{figure}[hbt!]
     \centering  
     \includegraphics[width=0.5\textwidth]{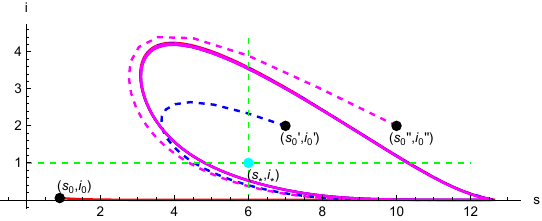}
       \caption{\footnotesize $(s,i)$-parametric plot with three different initial conditions illustrate the existence of stable limit cycles in blue, red, and magenta which collide around $(\xi,1)$ when $\xi=6$. The apparent intersection between the paths is due to the fact that we are projecting three dimensional paths in two dimensions.}
        \label{fig:1c}
    \end{figure}

\section{Application II: SIR{n}S models}\label{sec:sirns}

In this section, we discuss a SIRnS epidemic model for the evolution of susceptibles, infectives and recovered with several stages $s,i,r_1,...,r_n$, among a population that increases by a constant immigration rate $\lambda$, and decreases via natural deaths only. For the purposes of this paper, the model is somewhat complementary to the Capasso-Ruan-Wang models of Section \ref{s:SIRSgT} in the fact that here instabilities and bifurcations occur not due to different nonlinearities but rather just by the structural complexity of the network itself. In particular, we assume throughout mass action kinetics. The model is as follows:
  \begin{equation}\label{SIRn}
  \begin{cases} \dot{s}=\lambda - \mu s -\beta s i+\gamma_n r_n\\
i'=i (\beta s  + \sum_{k=1}^n \beta_k   r_k  -\gamma  -\mu )
\\\dot{r}_1=\gamma i - r_1(\beta_1    i+ \gamma_1   + \mu)\\ \vdots \\
\dot{r}_{k}=\gamma_{k-1} r_{k-1} - r_{k}(\beta_{k}    i + \gamma_k + \mu), k=2,...,n
\end{cases}.
\end{equation}

 \paragraph{Historical account.} The model \eqref{SIRn} above is motivated by the work of Herbert Hethcote and co-authors \cite{HethSVDD}, where the case of \eqref{SIRn} with $\beta=\beta_1=...\beta_n=\lambda=\mu_s=...\mu_{r_n}=0$ is considered.
Thus, they studied precisely a SIRnS model without demography, in which the loss of immunity is \emph{Erlang distributed}. Variations of the model with $\beta_k=0,  k =1, ...,n$ have also been studied in the context of enzyme kinetics, under the name of Glass-Kauffman model \cite{swick1985some,zaslavsky1995stability}. One interesting aspect of the model\cite{HethSVDD} is the occurrence of Hopf bifurcations, for $n\ge 3$. This is unusual for models which have a unique rational endemic point.
Subsequently, Adreu-Villaroig, González-Parra, and Villanueva \cite{Carlos} proposed introducing the parameters $\beta_k \ge 0$, which account for a gradual increase of infectiousness during the \emph{Erlang stages of transition} from $i$ to $s$;  this renders the problem more challenging: for example, the uniqueness of the endemic equilibrium is lost.

Epidemiological considerations thus advocate imposing the following constraints on the constants $\beta_k$:
\begin{equation}\label{eq:betas}
\beta>\beta_n>\beta_{n-1}>...>\beta_1.
\end{equation}
There is no known virus for which \eqref{eq:betas} doesn't hold. Nevertheless, the case of such `perverse viruses' deserves to be studied as well
Note that when $n=1$ and $\beta_1=0$, the model is a subcase of the Capasso-Ruan-Wang family, Section \ref{s:SIRSgT}.  
\

\subsection{General results}

For $n>3$, the system admits both zero-eigenvalue bifurcations and periodic orbits, even under the constraint \eqref{eq:betas}. We prove both cases in parallel via different methods: zero-eigenvalue bifurcations are achieved by a parametrization of the Jacobian at endemic equilibria, while periodic orbits are a consequence of the previous results by Hethcote and co-authors \cite{HethSVDD} and our inheritance Theorem \ref{thm:MEinheritance}.

\paragraph{Zero-eigenvalue bifurcations.}

We state and prove the following theorem.

\begin{theorem}\label{thm:sirnsn}
For any $n$, system \eqref{SIRn} admits endemic equilibria $(\bar{s},\bar{i},\bar{r}_1,...,\bar{r}_n)>0$ with a singular Jacobian. Moreover, for $n>3$, the system admits endemic equilibria with singular Jacobian even under the constraint \eqref{eq:betas}.
\end{theorem}

We postpone the lengthy proof of Theorem \ref{thm:sirnsn} in Appendix \ref{sec:thmproof}.

Theorem \ref{thm:sirnsn} guarantees the occurrence of zero-eigenvalue bifurcations for SIRnS systems with $n$ large enough. In turn, for $n=1$, constraint \eqref{eq:betas} excludes zero-eigenvalue bifurcations, as stated by the following Lemma.

\begin{lemma}\label{lemma:sir1s}
Let $n=1$. Then, for a zero-eigenvalue bifurcation, it is necessary that
\begin{equation}
\beta<\beta_1.
\end{equation}
\end{lemma}
The proof of Lemma \ref{lemma:sir1s} requires arguments from the proof of Theorem \ref{thm:sirnsn} and has been thus postponed as well in Appendix \ref{sec:thmproof}. We underline that Theorem \ref{thm:sirnsn} and Lemma \ref{lemma:sir1s} are not conclusive with respect of the cases $n=2,3$.\\

\paragraph{Hopf bifurcations.}

\begin{theorem}\label{thm:hopfsirns}
For $n\ge 3$, system \eqref{SIRn} admits  periodic orbits. Moreover, the system admits periodic orbits even under the constraints \eqref{eq:betas}.
\end{theorem}
The proof is based on inheritance results by Banaji and coauthors, Theorem \ref{thm:banaji} \cite{banaji2023bifurcation,banaji23split}, plus our own adaptation, Theorem \ref{thm:MEinheritance}. Such results are perturbation results, meaning that we prove the existence of periodic solutions only for a small enough choice of all the $\beta_k$. The persistence of such oscillatory behavior for large choices of $\beta_k$ must be further investigated.

\proof
The closed system with $\lambda=\mu=0$, and no re-infections $\beta_k=0$ for all $k$, admits nondegenerate Hopf bifurcations for $n\ge3$, as proved in \cite{HethSVDD}. Thus, we may consider at first the addition of nonzero $\beta_n$. The associated reaction
$$R_n+I\quad\underset{\beta_nir_n}\rightarrow \quad 2I$$
adds linearly dependent stoichiometry to the set of reactions
$$I\quad\rightarrow \quad R_1 \quad\rightarrow \quad R_2 \quad\rightarrow \quad ... \quad\rightarrow \quad R_n,$$
which are already present in the closed system \cite{HethSVDD}. Indeed:
the rank of the $n\times n$ stoichiometric submatrix matrix
\begin{equation}
\begin{pmatrix}
    -1 & 0 & 0 &...& 1\\
    1 & -1 & 0 & ...& 0\\
    0 & 1 & -1 & ... & 0\\
    & & & \vdots &\\
    0 & 0& 0& ... & -1
\end{pmatrix}
\end{equation}
is clearly $n-1$. Thus, the system with such an addition preserves the capacity for Hopf bifurcations, via the already mentioned Theorem \ref{thm:banaji}. We can then proceed by adding one by one reactions
$$R_{n-1}+I\rightarrow2I, \quad R_{n-2}+I\rightarrow2I, \quad ... \quad  R_{1}+I\rightarrow 2I.$$

It is crucial to follow this order: since Theorem \ref{thm:banaji} is a perturbation result, the result holds for $\beta_k$ in a right neighborhood of zero $\beta_k\in (0,\varepsilon_{\beta_k})$. This way, we are guaranteed that the constraint \eqref{eq:betas} can be enforced, simply by choosing $\varepsilon_{\beta_k}<\beta_{k+1}$.

Once all the infection-reactions have been added, we invoke our Theorem \ref{thm:MEinheritance} to conclude the result for general - and in principle independent - death rates and birth rate $\lambda$. To do so, we must only note that system \eqref{SIRn} is a travel-through system as defined in \ref{def:cyclic}. Again, because the inheritance results are perturbation results, we obtain the existence of periodic orbit for $\mu_k\in(0,\varepsilon_{\mu_k})$, and by taking the intersection of all such sets, we can consider $\mu_k\equiv\mu$ for all $k$ as in \eqref{SIRn}.

\endproof

\section{Conclusions and discussion}\label{s:Jin}
Our paper was motivated by the desire to investigate possible interactions between CRNT and ME. Our paper succeeded in finding several connections:
\begin{enumerate}
\item We imported from CRNT the idea \cite{Vas} to split the search for bifurcation points into two steps: the first, reflecting only the network structure,  treats the reaction rates $\mathbf{r}(x)$ and their derivatives as free variables. The resulting `symbolic bifurcation problem' is naturally easier to solve. In this direction, Theorem \ref{lem:purelyimaginary} states that essentially all epidemic models admit `symbolic Hopf bifurcations'. To connect to the classic parametric problems, one assumes then a specific parametric form for $\mathbf{r}(x)$, often mass-action kinetics, and solves the resulting algebraic problem of finding parameter values which fit the result of the first step. Here we may conclude that
all epidemic models endowed with Michaelis-Menten rates admit an endemic equilibrium with a purely-imaginary-eigenvalues Jacobian. It is only the restriction to mass action that may prevent the occurrence of Hopf bifurcations.
\item The previous ideas seemed hand-made for the class of Capasso-Ruan-Wang models, Section \ref{s:SIRSgT}, for which we succeeded in providing a unified result, Theorem \ref{thm:necessarybif}, which had been partially achieved separately in many papers for each specification of the model. Of course, the original papers had obtained sharper results by their case-to-case analysis, as that was the very focus of such contributions.
\item  We further demonstrated the utility of CRNT methods in the recent (2n+d) SIRnS system  \cite{Carlos}. For a purely mass-action case, we have obtained a parametrization of the equilibria and proved subsequently the existence of endemic equilibria undergoing a zero-eigenvalue bifurcation, for every $n$. As an aside, we note that while rational parametrizations have been employed in ME, they have again been derived case by case. In contrast, the CRNT parametrization reveals from the start the connection to the null space of the stoichiometric matrix $\Gamma$.
\item Finally, we have built on Banaji's inheritance result \cite{banaji23split}, and extended it in Theorem \ref{thm:MEinheritance} to prove that Hopf bifurcations in closed `travel-through' models imply periodic orbits also for the open model with demographics. This result allowed us to prove the existence of Hopf bifurcations for the epidemic model in the SIRnS system \cite{Carlos} with $n\geq 3$, simply by relying on the previous result for a simpler model \cite{HethSVDD}.
\end{enumerate}

 In conclusion, our results support the utility of
adopting CRNT methods in the study of ME models.
For interactions in the opposite direction, future work can use methods from ME to study the loss of stability of boundary equilibria of biochemical networks. As well known, the general structure available in CRNT models has allowed very fruitful interactions with researchers in real algebraic geometry, giving rise to the emerging field of algebraic biology. We hope that such interactions will also be possible for ME models.

\textbf{ Acknowledgements.}  We thank Carlos Andreu, Gilberto Gonzalez-Parra, and Rafael Villanueva for sharing with us their beautiful recent model addressing the issue of "immunity in influenza epidemiological dynamics", as well as for providing the crucial reference to the work of Hethcote,  Stech, and  Van Den Driessche \cite{HethSVDD}. We also thank Murad Banaji, Balasz Boros, Dan Goreac, Daniel Lichtblau, Matteo Levi, Mattia Sensi, Carmen Rocsoreanu, Anne Shiu and Janos Toth for useful advice. Nicola Vassena has been supported by the DFG (German Research
Society), project n. 512355535.



\newpage
\appendix
\section{Appendix: Proof of Theorem \ref{thm:sirnsn} and Lemma \ref{lemma:sir1s}}\label{sec:thmproof}

We present the proof of Theorem \ref{thm:sirnsn} divided into four lemmas. The first presents a parametrization of the vector of the equilibrium flux cone, i.e. the reaction rates evaluated at any endemic equilibrium. Let $$d_s,d_i,d_1...,d_{n},b_1,...,b_n,c>0$$ be $2n+3$ positive parameters.

\begin{lemma}[Parametrization of the flux cone]\label{lemma:z1}
     At an endemic equilibrium $(\bar{s},\bar{i},\bar{r}_1,...,\bar{r}_n)$, the vector of the reaction rates can be parameterized as follows:
\begin{equation}
 \begin{pmatrix}
           \lambda\\
           \beta \bar{s}\bar{i}\\
    \beta_1 \bar{i} \bar{r_1}\\
    \vdots\\
    \beta_n \bar{i} \bar{r_n}\\
    \gamma \bar{i}\\
    \gamma_1 \bar{r}_1\\
    \vdots\\
    \gamma_k \bar{r}_k\\
    \vdots\\
    \gamma_n \bar{r}_n\\
    \mu_s \bar{s}\\
    \mu_i \bar{i}\\
    \mu_1 \bar{r_1}\\
    \vdots\\
    \mu_n \bar{r_n}\\
\end{pmatrix} =\begin{pmatrix}
d_s+d_i+\sum_{k=1}^n d_k\\
c+d_i+\sum_{k=1}^n d_k\\
b_1\\
\vdots\\
b_{n}\\
c+\sum_{j=1}^n(b_{j}+d_{j})\\
c+\sum_{j=2}^n(b_{j}+d_{j})\\
\vdots\\
c+\sum_{j=k+1}^n (b_{j}+d_{j})\\
\vdots\\
c\\
d_s\\
d_i\\
d_1\\
\vdots\\
d_{n}\\
\end{pmatrix}
    \end{equation}
\end{lemma}
\proof
The proof is by direct computation. We start from the last equation of $r'_n$:
\begin{equation}0=\gamma_{n-1}\bar{r}_{n-1}-\beta_n \bar{i} \bar{r}_n-\gamma_n\bar{r}_n -\mu_n \bar{r}_n.
\end{equation}
We define the free parameters $b_n:=\beta_n \bar{i} r_n,$, $c:=\gamma_n\bar{r}_n$, $d_n:=\mu_n \bar{r}_n$, which leads to
\begin{equation}
  \gamma_{n-1}\bar{r}_{n-1}=c+b_n+d_n.
\end{equation}
Since $\gamma_{n-1}\bar{r}_{n-1}$ appears also in the equation for $r'_{n-1}$ we can analogously define the free parameters $b_{n-1}:=\beta_{n-1}\bar{i} \bar{r}_{n-1}$, $d_{n-1}:=\mu_{n-1} \bar{r}_{n-1}$ and re-iterate the argument obtaining
\begin{equation}\gamma_{n-2}r_{n-2}=b_{n-1}-\gamma_{n-1}\bar{r}_{n-1}-d_{n-1}=b_{n-1}-c+b_n+d_n-d_{n-1}.
\end{equation}
In general, we can define $b_k:=\beta_k i r_k$, $d_k:=\mu_k r_k$ and via a straightforward iteration
\begin{equation}\gamma_{k}r_{k}=c+\sum_{j=k+1}^n(b_j+d_j)\quad\text{and}\quad \gamma i = c+\sum_{j=1}^n+(b_j+d_j)
\end{equation}
Now, we only need to satisfy equilibria equations for $s$ and $i$. Define further $d_s:=\mu_s \bar{s}, d_i:=\mu_i \bar{i}$. From the $i$-equation we get
\begin{equation}
\begin{split}
0=&\beta \bar{s} \bar{i} + \sum_{k=1}^n b_k - \gamma \bar{i} - d_i \\
=&\beta \bar{s} \bar{i} + \sum_{k=1}^n b_k - c - \sum_{k=1}^n(b_k+d_k) - d_i\\
=&\beta \bar{s} \bar{i} - c - d_i - \sum_{k=1}^n d_k,
\end{split}
\end{equation}
from which we get the equality
$$\beta \bar{s} \bar{i} = c + d_i + \sum_{k=1}^n d_k.$$
Finally, from the equation for $s$ we get
\begin{equation}
    \begin{split}
        0=& \lambda-\beta \bar{s} \bar{i} + c - d_s\\
        =& \lambda -c - d_i - \sum_{k=1}^n d_k +c -d_s\\
        =& \lambda - d_s - d_i - \sum_{k=1}^n d_k,
    \end{split}
\end{equation}
obtaining
$$\lambda = d_s + d_i + \sum_{k=1}^n d_k.$$
\endproof
The second lemma parametrizes the Jacobian matrix evaluated at endemic equilibria.
\begin{lemma}[Parametrization of the Jacobian at endemic equilibria]\label{lemma:sirn2}
  At an endemic equilibrium $\bar{x}:=(\bar{s},\bar{i},\bar{r}_1,...,\bar{r}_n)>0$ the Jacobian $G(\bar{x})$ can be parametrized as follows:
{\footnotesize
\begin{equation}\label{eq:sirnparjac}
G(b_1,...,b_n,c,d_s,d_i,d_1...,d_{n};\bar{s},\bar{i},\bar{r}_1,...,\bar{r}_n)= J(b_1,...,b_n,c,d_s,d_i,d_1...,d_{n})\operatorname{diag}\bigg(\frac{1}{\bar{s}},\frac{1}{\bar{i}},\frac{1}{\bar{r}_1},...,\frac{1}{\bar{r}_n}\bigg),\end{equation}}
where $J(b_1,...,b_n,c,d_s,d_i,d_1...,d_{n})$ is the following matrix:

{\scriptsize
\begin{equation*}
\begin{pmatrix}
-c-d_s-d_i-\sum_{k=1}^n d_k & -c-d_i-\sum_{k=1}^n d_k & 0 & 0 & ... & c\\
c+d_i+\sum_{k=1}^n d_k& 0 & b_1 & b_2 & ... & b_n\\
0 & c+d_1+\sum_{j=2}^n(b_{j}+d_{j}) & -c-\sum_{j=1}^n(b_{j}+d_{j}) & 0 & ... & 0\\
0 & -b_2 & c+\sum_{j=2}^n(b_{j}+d_{j}) & -c-\sum_{j=2}^n(b_{j}+d_{j}) & ... & 0\\
0 & -b_3 & 0 &  c+\sum_{j=3}^n(b_{j}+d_{j})& ... & 0\\
& & & & \vdots & \\
0 & -b_n & 0 & 0 &... & -c-b_n-d_n\\
\end{pmatrix}.
\end{equation*}}
\end{lemma}

\proof
First, we derive the general form of the Jacobian $G$ by explicit differentiation.
\begin{equation*}
\begin{split}
G=\begin{pmatrix}
\beta \bar{i}-\mu_s & \beta \bar{s} & 0 & 0 & ... & \gamma_n\\
\beta \bar{i} & 0 & \beta_1 \bar{i} & \beta_2 \bar{i} & ... & \beta_n \bar{i}\\
0 & \gamma - \beta_1 \bar{r}_1 & -\gamma_1 - \beta_1 \bar{i}-\mu_1 & 0 & ... & 0\\
0 & -\beta_2 \bar{r}_2 & \gamma_1 & -\gamma_2 - \beta_2 \bar{i}-\mu_2 & ... & 0\\
0 & -\beta_3 \bar{r}_3 & 0 & \gamma_2& ... & 0\\
& & & & \vdots & \\
0 & -\beta_n \bar{r}_n & 0 & 0 & ... & -\gamma_n-\beta_n \bar{i}-\mu_n\\
\end{pmatrix}
\end{split}
\end{equation*}
Second, we just consider
$$G=G \operatorname{Id}= G \;\operatorname{diag}(\bar{s},\bar{i},\bar{r}_n,...,\bar{r}_n) \; \operatorname{diag}\bigg(\frac{1}{\bar{s}},\frac{1}{\bar{i}},\frac{1}{\bar{r}_1},...,\frac{1}{\bar{r}_n}\bigg),$$
where we note that the entries of
$$J=G \operatorname{diag}(\bar{s},\bar{i},\bar{r}_n,...,\bar{r}_n)$$
are simply a sum of reaction rates evaluated at the equilibrium. We can then use the parametrization of Lemma 1 to write $J$ in the above form.
\endproof

The third lemma states that there are choices of parameters for which the Jacobian has a sign-change of the determinant, hereby implying a zero-eigenvalue bifurcation.
\begin{lemma}[Change of equilibria stability]\label{lemma:z3}
The determinant of the Jacobian at an endemic equilibrium can be expressed as:
$$\operatorname{det}G=B-Ad_s,$$
where the terms $A,B$ do not depend on $d_s$ and can be chosen of sign $(-1)^n$. Thus, there are parameter values for which the determinant changes sign.
\end{lemma}
\proof
Consider the Jacobian parametrization \eqref{eq:sirnparjac} of Lemma \ref{lemma:sirn2} above. Proving the statement for $J$ alone is sufficient, since the
$$\operatorname{det}G=\operatorname{det}J\; (\bar{s}\bar{i}\bar{r}_1...\bar{r}_n)^{-1}.$$
Firstly note that that for $j > 3$,
$$J_{jj}=-J_{j(j-1)}.$$
We operate then the following Gauss-type simplification:
\begin{enumerate}
\item $J^1$ is the matrix obtained from $J$ by summing the $n^{th}$ column of $J$ to its $(n-1)^{th}$ column;
\item $J^2$ is the matrix obtained by $J^1$ by summing the $(n-1)^{th}$ column of $J^1$ to its $(n-2)^{th}$ column
\item ...
\item $J^k$, $k=3,...,n-2$ is the matrix obtained by $J^{k-1}$ by summing the $n-k+1$ column to its $(n-k)^{th}$ column.
\end{enumerate}
To simplify the notation, define $$\tilde{J}:=J^{n-2}$$
Clearly, $\operatorname{det}(\tilde{J})=\operatorname{det}(J)$. However, $\tilde{J}$ is now in the convenient form:
{\scriptsize$$\begin{pmatrix}
-c-d_s-d_i-\sum_{k=1}^n d_k & -d_i-\sum_{k=1}^n d_k & c & c & ... & c\\
c+d_i+\sum_{k=1}^n d_k& \sum_{k=1}^n b_k & \sum_{k=1}^n b_k & \sum_{k=2}^n b_k & ... & b_n\\
0 & -b_1 & -c-\sum_{j=1}^n(b_{j}+d_{j}) & 0 & ... & 0\\
0 & -b_2 & 0 & -c-\sum_{j=2}^n(b_{j}+d_{j}) & ... & 0\\
0 & -b_3 & 0 &  0& ... & 0\\
& & & & \vdots & \\
0 & -b_n & 0 & 0 &... & -c-b_n-d_n\\
\end{pmatrix}.$$}
Note that $d_s$ appears only in the entry $\tilde{J}_{11}$. Now, define the four blocks
\begin{equation}
    \begin{split}
        A&:=\begin{pmatrix}
   -c-d_s-d_i-\sum_{k=1}^n d_k & -d_i-\sum_{k=1}^n d_k\\
  c+d_i+\sum_{k=1}^n d_k& \sum_{k=1}^n b_k\\
  \end{pmatrix}\\
  B&:=\begin{pmatrix}
   c & c & ... & c\\
   \sum_{k=1}^n & \sum_{k=2}^n & ... & b_n\\
\end{pmatrix}\\
C&:=\begin{pmatrix}
    0 & -b_1\\
    0 & -b_2\\
    \vdots& \vdots\\
    0 & -b_n\\
\end{pmatrix}\\
D&:=\begin{pmatrix}
-c-\sum_{j=1}^n(b_{j}+d_{j}) & 0 & ... & 0\\
0 & -c-\sum_{j=2}^n(b_{j}+d_{j}) & ... & 0\\
& & \vdots &\\
0 & 0 & 0 & -c -b_n-d_n
\end{pmatrix}
    \end{split}
\end{equation}

From the fact that $$\tilde{J}=\begin{pmatrix}
A & B\\
C & D\\
\end{pmatrix},
$$
and that $D$ is always invertible, we get
$$\operatorname{sign} \operatorname{det}(\tilde{J})=\operatorname{sign} \operatorname{det}(D) \operatorname{det}(A-BD^{-1}C)=(-1)^n \operatorname{sign} \operatorname{det}(A-BD^{-1}C).$$
Again, proving the statement for $(A-BD^{-1}C)$ alone is sufficient to conclude. Ugly-looking, but straightforward computation yields:
{\footnotesize
\begin{equation*}
    \begin{split}
    BD^{-1}C&=\begin{pmatrix}
   c & c & ... & c\\
   \sum_{k=1}^n b_k & \sum_{k=2}^n b_k& ... & b_n\\
\end{pmatrix} \begin{pmatrix}
\frac{1}{-c-\sum_{j=1}^n(b_{j}+d_{j})} & 0 & ... & 0\\
0 & \frac{1}{-c-\sum_{j=2}^n(b_{j}+d_{j})} & ... & 0\\
& & \vdots &\\
0 & 0 & 0 & \frac{1}{-c -b_n-d_n}
\end{pmatrix}
\begin{pmatrix}
    0 & -b_1\\
    0 & -b_2\\
    \vdots& \vdots\\
    0 & -b_n\\
\end{pmatrix}\\
&=
\begin{pmatrix}
\frac{c}{-c-\sum_{j=1}^n(b_{j}+d_{j})} & \frac{c}{-c-\sum_{j=2}^n(b_{j}+d_{j})} & ... & \frac{c}{-c-b_{n}-d_{n}}\\
\frac{\sum_{k=1}^n b_k}{-c-\sum_{j=1}^n(b_{j}+d_{j})} & \frac{\sum_{k=2}^n b_k}{-c-\sum_{j=2}^n(b_{j}+d_{j})} & ... &\frac{ b_n}{-c-b_n-d_n}\\
\end{pmatrix}\begin{pmatrix}
    0 & -b_1\\
    0 & -b_2\\
    \vdots& \vdots\\
    0 & -b_n\\
\end{pmatrix}\\
&=\begin{pmatrix}
0 & \sum_{k=1}^n\frac{b_k c}{c+\sum_{j=k}^n(b_{j}+d_{j})}\\
0 & \sum_{k=1}^n \frac{b_k\sum_{j=k}^n b_j}{c+\sum_{j=k}^n(b_{j}+d_{j})}\\
\end{pmatrix}.
\end{split}
\end{equation*}}
For simplicity, set \begin{equation}\label{eq:alpha1}
\alpha_1:=\sum_{k=1}^n\frac{b_k c}{c+\sum_{j=k}^n(b_{j}+d_{j})},
\end{equation}and
\begin{equation}\label{eq:alpha2}
\alpha_2:=\sum_{k=1}^n \frac{b_k\sum_{j=k}^n b_j}{c+\sum_{j=k}^n(b_{j}+d_{j})}.
\end{equation}
Note that
$$0<\alpha_1,\alpha_2<\sum_{k=1}^n b_k$$
For the positive choices of
\begin{equation}\label{eq:deltaalpha34}
    \begin{cases}
\delta:=c+d_1+\sum_{k=1}^n d_k>0\\
    \alpha_3:= d_i + \sum_{k=1}^n d_k +\alpha_1>0\\
\alpha_4:=\sum_{k=1}^n b_k-\alpha_2>0\\
    \end{cases},
\end{equation}
we get:
\begin{equation}\label{eq:detjac}
    \begin{split}
        \Delta:=\operatorname{det}(A-BD^{-1}C)&=\operatorname{det}\begin{pmatrix}
   -c-d_s-d_i-\sum_{k=1}^n d_k & -d_i-\sum_{k=1}^n d_k-\alpha_1\\
  c+d_i+\sum_{k=1}^n d_k& \sum_{k=1}^n b_k-\alpha_2\\
  \end{pmatrix}\\
  &=\operatorname{det}\begin{pmatrix}
  -d_s-\delta & -\alpha_3\\
  \delta & \alpha_4,\\
  \end{pmatrix}
\\
&=-d_s\alpha_4+\delta (\alpha_3-\alpha_4 ).
    \end{split}
\end{equation}
We need now to show only that $\alpha_3-\alpha_4$ can be chosen positive. We compute:
\begin{equation}
\begin{split}
\alpha_3-\alpha_4&= d_i + \sum_{k=1}^n d_k +\alpha_1-\sum_{k=1}^n b_k+\alpha_2\\
&=d_i + \sum_{k=1}^n d_k + \sum_{k=1}^n \bigg(b_k \bigg(\frac{c+\sum_{j=k}^n b_j}{c+\sum_{j=k}^n (b_j+d_j)}-1 \bigg)\bigg)\\
&=d_i + \sum_{k=1}^n d_k + \sum_{k=1}^n \bigg(b_k \bigg(\frac{c+\sum_{j=k}^n b_j+d_j-d_j}{c+\sum_{j=k}^n (b_j+d_j)}-1 \bigg)\bigg)\\
&=d_i+\sum_{k=1}^n d_k + \sum_{k=1}^n \bigg(b_k \bigg(-\frac{\sum_{j=k}^nd_j}{c+\sum_{j=k}^n (b_j+d_j)} \bigg)\bigg)\\
&=d_i+\sum_{k=1}^n d_k-\sum_{k=1}^n \bigg(d_k \bigg(\frac{\sum_{j=1}^kb_j}{c+\sum_{j=k}^n (b_j+d_j)} \bigg)\bigg)\\
&=d_i+\sum_{k=1}^n\bigg(d_k\bigg(1- \frac{\sum_{j=1}^kb_j}{c+\sum_{j=k}^n (b_j+d_j)} \bigg)\bigg).
\end{split}
\end{equation}
Clearly, we may choose the parameters $b_k$, $k=1,...,n$, small enough so that $$\alpha_3-\alpha_4>0.$$
\endproof

\begin{lemma}\label{lemma:z4}
For $n> 3$ there exists a positive choice of parameters satisfying constraint \eqref{eq:betas} for which there is an equilibrium with singular Jacobian.
\end{lemma}
\proof
Based on the parametrization of Lemma \ref{lemma:sirn2}, we will find two explicit positive choices $\bar{\mathbf{p}}_1>0$ and $\bar{\mathbf{p}}_2>0$ of parameters such that, for any $n>3$, the equilibrium Jacobian \eqref{eq:sirnparjac} has determinant of sign $(-1)^n$ for $\bar{\mathbf{p}}_1$ and of sign $(-1)^{n-1}$ for $\bar{\mathbf{p}}_2$. The intermediate value theorem directly implies the existence of a choice $\bar{\mathbf{p}}^*$ of parameters where the equilibrium has a zero-eigenvalue determinant. We argue as follows. Firstly fix
\begin{equation}\label{eq:parsir6s}
\gamma=\gamma_1=...\gamma_n=\mu=1,\quad\text{ as well as}\quad \bar{s}=\bar{i}=\bar{r}_1=...=\bar{r}_n=1.
\end{equation}
As a consequence of \eqref{eq:parsir6s}, we have now also implicitly fixed the related parameters defined in Lemma \ref{lemma:z1}:
$$d_s=d_i=d_1=...=d_n=c=1.$$
and, crucially, we also get that
$$\beta=\beta\bar{s}\bar{i}=d_i+d_1+...+d_n+c=n+2.$$
We keep $\beta_1,...,\beta_n$ as free parameters, and consequently
$$b_k=\beta_k,\quad\text{for $k=1,...,n$}.$$
If we show that there exists a positive choice of parameters satisfying
\begin{equation}\label{eq:betabcons}
\beta_1=\beta_2=...=\beta_n=:\bar{b}\le \beta=n+2,
\end{equation}
we then easily conclude for the statement of Lemma \ref{lemma:z4}. In fact, simply define then
$$\beta_k=\bar{b}-\frac{\varepsilon}{k},$$
and obtain the statement of Lemma \ref{lemma:z4} for $\varepsilon$ small enough.

We then proceed assuming constraint \eqref{eq:betabcons}. We first recall, via the proof of Lemma \ref{lemma:z3}, that the sign of the determinant uniquely depends on the expression \eqref{eq:detjac}
$$\Delta=-d_s\alpha_4+\delta(\alpha_3-\alpha_4).$$
In particular, if $\Delta>0$ then the determinant of the Jacobian is of sign $(-1)^n$ whereas if $\Delta<0$ the determinant is of sign $(-1)^{n-1}$. We again underline that for $\bar{b}>0$ small enough,
$$\operatorname{sign}\Delta\approx \operatorname{sign}\delta\alpha_3>0,$$
and thus we always have a positive choice of parameters such that the Jacobian determinant is of sign $(-1)^{n}$.

In turn, we evaluate $\Delta$ at the values  \eqref{eq:parsir6s} and \eqref{eq:betabcons} with $\bar{b}=n+2$. Note that, at such values:
\begin{equation}
\begin{cases}
\delta=n+2;\\
\alpha_4=n(n+2)-\sum_{k=1}^n\frac{(n+2)(n+2)(n-k+1)}{1+(n-k+1)(n+3)};\\
\alpha_3-\alpha_4=1+n-\sum_{k=1}^{n}\frac{k(n+2)}{1+(n-k+1)(n+3)}.
\end{cases}
\end{equation}
We then compute $\Delta=-\alpha_4+\delta(\alpha_3-\alpha_4)$:
\begin{equation}
\begin{split}
\Delta(n)&=(n+2)\bigg(-n+\sum_{k=1}^n\frac{(n+2)(n-k+1)}{1+(n-k+1)(n+3)}+ 1+n-\sum_{k=1}^{n}\frac{k}{1+(n-k+1)(n+3)}   \bigg)\\
&=(n+2)\bigg(1+(n+2)\sum_{k=1}^n\frac{n-2k+1}{1+(n-k+1)(n+3)}\bigg),
\end{split}
\end{equation}
whose sign depends only on $$\varsigma(n):=1+(n+2)\sum_{k=1}^n\dfrac{n-2k+1}{1+(n-k+1)(n+3)}.$$
A straightforward computation shows that $\varsigma(n)<0$ if $n=4$. Since $\varsigma(n)$ is monotone decreasing, we get that $\varsigma(n)<0$ for all $n>3$. This concludes the proof of Lemma \ref{lemma:z4}.
\endproof

The theorem is proven, as a direct consequence of Lemmas \ref{lemma:z3} and \ref{lemma:z4}.

We conclude this Section with the proof of Lemma \ref{lemma:sir1s}.

\proof[Proof of Lemma \ref{lemma:sir1s}]
Firstly, we recall the values $\delta,\alpha_1,\alpha_2,\alpha_3,\alpha_4$ from \eqref{eq:alpha1}, \eqref{eq:alpha2}, \eqref{eq:deltaalpha34}. We note the following equalities:
\begin{equation}
\beta=\frac{\delta}{\bar{s}\bar{i}}=\frac{\mu^2 \delta}{d_s d_i};\quad \text{and}\quad
\beta_1=\frac{b_1}{\bar{r}_1\bar{i}}=\frac{\mu^2 b_1}{d_1 d_i},
\end{equation}
from which we first derive that
\begin{equation}\label{eq:betabeta}
    \beta \beta_1^{-1}=\frac{\delta d_1}{d_s b_1}.
\end{equation}

Note that for $n=1$,
$$\alpha_3-\alpha_4=d_i+d_1\bigg(1-\frac{b_1}{c+b_1+d_1}\bigg)>0,$$
and
$$\alpha_4=b_1\bigg(1-\frac{b_1}{c+b_1+d_1}\bigg).$$
In particular, a zero-eigenvalue is achieved if and only if
$$d^*_s:=\delta \frac{(\alpha_3-\alpha_4)}{\alpha_4}.$$

We then estimate $\delta/d_s$ at the bifurcation point $d^*_s$.
$$\frac{\delta}{d^*_s}=\frac{\alpha_4}{\alpha_3-\alpha_4}=\frac{b_1\bigg(1-\frac{b_1}{c+b_1+d_1}\bigg)}{d_i+d_1\bigg(1-\frac{b_1}{c+b_1+d_1}\bigg)}<\frac{b_1\bigg(1-\frac{b_1}{c+b_1+d_1}\bigg)}{d_1\bigg(1-\frac{b_1}{c+b_1+d_1}\bigg)}=\frac{b_1}{d_1}.$$

Substituting in \eqref{eq:betabeta} we get $$\beta< \beta_1$$

\endproof

\bibliography{references.bib}

\providecommand{\bysame}{\leavevmode\hbox to3em{\hrulefill}\thinspace}
\providecommand{\MR}{\relax\ifhmode\unskip\space\fi MR }
\providecommand{\MRhref}[2]{%
  \href{http://www.ams.org/mathscinet-getitem?mr=#1}{#2}
}
\providecommand{\href}[2]{#2}
\begin{thebibliography}{HSVDD81}

\bibitem[ADLS07]{Ang07}
David Angeli, Patrick De~Leenheer, and Eduardo~D Sontag, \emph{A {P}etri net
  approach to the study of persistence in chemical reaction networks},
  Mathematical biosciences \textbf{210} (2007), no.~2, 598--618.

\bibitem[AGPV24]{Carlos}
Carlos Andreu, Gilberto Gonzalez-Parra, and Rafael Villanueva, \emph{In
  preparation}.

\bibitem[Ang09]{Ang}
David Angeli, \emph{A tutorial on chemical reaction networks dynamics}, 2009
  European Control Conference (ECC), IEEE, 2009, pp.~649--657.

\bibitem[ARAS23]{al2023structural}
M~Ali Al-Radhawi, David Angeli, and Eduardo Sontag, \emph{On structural
  contraction of biological interaction networks}, arXiv preprint
  arXiv:2307.13678 (2023).

\bibitem[Ban23]{banaji23split}
Murad Banaji, \emph{Splitting reactions preserves nondegenerate behaviors in
  chemical reaction networks}, SIAM Journal on Applied Mathematics \textbf{83}
  (2023), no.~2, 748--769.

\bibitem[BB23]{BaBo23}
Murad Banaji and Bal{\'a}zs Boros, \emph{The smallest bimolecular mass action
  reaction networks admitting andronov--{H}opf bifurcation}, Nonlinearity
  \textbf{36} (2023), no.~2, 1398.

\bibitem[BBH23]{banaji2023bifurcation}
Murad Banaji, Bal{\'a}zs Boros, and Josef Hofbauer, \emph{The inheritance of
  local bifurcations in mass action networks}, arXiv preprint arXiv:2312.12897
  (2023).

\bibitem[BC09]{banaji2009graph}
Murad Banaji and Gheorghe Craciun, \emph{Graph-theoretic approaches to
  injectivity and multiple equilibria in systems of interacting elements}.

\bibitem[BP16]{BaCa}
Murad Banaji and Casian Pantea, \emph{Some results on injectivity and
  multistationarity in chemical reaction networks}, SIAM Journal on Applied
  Dynamical Systems \textbf{15} (2016), no.~2, 807--869.

\bibitem[CDSS09]{CDSS}
Gheorghe Craciun, Alicia Dickenstein, Anne Shiu, and Bernd Sturmfels,
  \emph{Toric dynamical systems}, Journal of Symbolic Computation \textbf{44}
  (2009), no.~11, 1551--1565.

\bibitem[CF05]{CF05}
Gheorghe Craciun and Martin Feinberg, \emph{Multiple equilibria in complex
  chemical reaction networks: I. the injectivity property}, SIAM Journal on
  Applied Mathematics \textbf{65} (2005), no.~5, 1526--1546.

\bibitem[Cla88]{ClarkeSNA}
Bruce~L Clarke, \emph{Stoichiometric network analysis}, Cell biophysics
  \textbf{12} (1988), 237--253.

\bibitem[Cox20]{cox2020applications}
David~A Cox, \emph{Applications of polynomial systems}, vol. 134, American
  Mathematical Soc., 2020.

\bibitem[CP19]{CoCa}
Carsten Conradi and Casian Pantea, \emph{Multistationarity in biochemical
  networks: results, analysis, and examples}, Algebraic and combinatorial
  computational biology, Elsevier, 2019, pp.~279--317.

\bibitem[CS78]{Capasso}
Vincenzo Capasso and Gabriella Serio, \emph{A generalization of the
  {K}ermack-{M}c{K}endrick deterministic epidemic model}, Mathematical
  biosciences \textbf{42} (1978), no.~1-2, 43--61.

\bibitem[DFS20]{degrand2020graphical}
Elisabeth Degrand, Fran{\c{c}}ois Fages, and Sylvain Soliman, \emph{Graphical
  conditions for rate independence in chemical reaction networks},
  International Conference on Computational Methods in Systems Biology,
  Springer, 2020, pp.~61--78.

\bibitem[DHM90]{Diek}
Odo Diekmann, Johan Andre~Peter Heesterbeek, and Johan~AJ Metz, \emph{On the
  definition and the computation of the basic reproduction ratio {R}0 in models
  for infectious diseases in heterogeneous populations}, Journal of
  mathematical biology \textbf{28} (1990), no.~4, 365--382.

\bibitem[Fei87]{Fei87}
Martin Feinberg, \emph{Chemical reaction network structure and the stability of
  complex isothermal reactors—{I}. {T}he deficiency zero and deficiency one
  theorems}, Chemical Engineering Science \textbf{42} (1987), no.~10,
  2229--2268.

\bibitem[FGS15]{fages2015inferring}
Fran{\c{c}}ois Fages, Steven Gay, and Sylvain Soliman, \emph{Inferring reaction
  systems from ordinary differential equations}, Theoretical Computer Science
  \textbf{599} (2015), 64--78.

\bibitem[Fie85]{fiedler1985index}
Bernold Fiedler, \emph{An index for global {H}opf bifurcation in parabolic
  systems.}, Journal f{\"u}r die reine und angewandte Mathematik \textbf{358}
  (1985), 1--36.

\bibitem[GH84]{GuHo84}
John Guckenheimer and Philip Holmes, \emph{Nonlinear oscillations, dynamical
  systems and bifurcations of vector fields}, Springer, 1984.

\bibitem[GK22]{Gupta}
RP~Gupta and Arun Kumar, \emph{Endemic bubble and multiple cusps generated by
  saturated treatment of an {SIR} model through {H}opf and {B}ogdanov--{T}kens
  bifurcations}, Mathematics and Computers in Simulation (2022).

\bibitem[GS23]{golubitsky2023dynamics}
Martin Golubitsky and Ian Stewart, \emph{Dynamics and bifurcation in networks:
  Theory and applications of coupled differential equations}, SIAM, 2023.

\bibitem[HCH10]{Haddad}
Wassim~M Haddad, VijaySekhar Chellaboina, and Qing Hui, \emph{Nonnegative and
  compartmental dynamical systems}, Princeton University Press, 2010.

\bibitem[Hil10]{HIll10}
Archibald~Vivian Hill, \emph{The possible effects of the aggregation of the
  molecules of haemoglobin on its dissociation curves}, The Journal of
  Physiology \textbf{40} (1910), 4--7.

\bibitem[HJ72]{horn1972general}
Fritz Horn and Roy Jackson, \emph{General mass-action kinetics}, Archive for
  rational mechanics and analysis \textbf{47} (1972), 81--116.

\bibitem[Hol65]{Holling65}
Crawford~Stanley Holling, \emph{The functional response of predators to prey
  density and its role in mimicry and population regulation}, The Memoirs of
  the Entomological Society of Canada \textbf{97} (1965), no.~S45, 5--60.

\bibitem[HSVDD81]{HethSVDD}
Herbert~W Hethcote, Harlan~W Stech, and Pauline Van Den~Driessche,
  \emph{Nonlinear oscillations in epidemic models}, SIAM Journal on Applied
  Mathematics \textbf{40} (1981), no.~1, 1--9.

\bibitem[LHL87]{LiuLevin}
Wei-min Liu, Herbert~W Hethcote, and Simon~A Levin, \emph{Dynamical behavior of
  epidemiological models with nonlinear incidence rates}, Journal of
  mathematical biology \textbf{25} (1987), 359--380.

\bibitem[LHRY19]{LuRuan19}
Min Lu, Jicai Huang, Shigui Ruan, and Pei Yu, \emph{Bifurcation analysis of an
  {SIRS} epidemic model with a generalized nonmonotone and saturated incidence
  rate}, Journal of differential equations \textbf{267} (2019), no.~3,
  1859--1898.

\bibitem[Lot20]{Lotka20}
Alfred~J Lotka, \emph{Analytical note on certain rhythmic relations in organic
  systems}, Proceedings of the National Academy of Sciences \textbf{6} (1920),
  no.~7, 410--415.

\bibitem[MM13]{MM13}
L.~Michaelis and M.~L. Menten, \emph{Die kinetik der invertinwirkung}, Biochem.
  Z. \textbf{49} (1913), 333--369.

\bibitem[MY20]{macauley2020case}
Matthew Macauley and Nora Youngs, \emph{The case for algebraic biology: from
  research to education}, Bulletin of Mathematical Biology \textbf{82} (2020),
  1--16.

\bibitem[PS05]{pachter2005algebraic}
Lior Pachter and Bernd Sturmfels, \emph{Algebraic statistics for computational
  biology}, vol.~13, Cambridge university press, 2005.

\bibitem[RBG22]{Roostaei}
Arash Roostaei, Hadi Barzegar, and Fakhteh Ghanbarnejad, \emph{Emergence of
  {H}opf bifurcation in an extended {SIR} dynamic}, Plos one \textbf{17}
  (2022), no.~10, e0276969.

\bibitem[SF12]{ShiFei12}
Guy Shinar and Martin Feinberg, \emph{Concordant chemical reaction networks},
  Mathematical Biosciences \textbf{240} (2012), no.~2, 92--113.

\bibitem[SH81]{sokol81}
W~Sokol and JA~Howell, \emph{Kinetics of phenol oxidation by washed cells},
  Biotechnology and Bioengineering \textbf{23} (1981), no.~9, 2039--2049.

\bibitem[Sol13]{soliman2013stronger}
Sylvain Soliman, \emph{A stronger necessary condition for the multistationarity
  of chemical reaction networks}, Bulletin of mathematical biology \textbf{75}
  (2013), 2289--2303.

\bibitem[Swi85]{swick1985some}
KE~Swick, \emph{Some reducible models of age dependent dynamics}, SIAM Journal
  on Applied Mathematics \textbf{45} (1985), no.~2, 256--267.

\bibitem[TF21]{torres2021symbolic}
Ang{\'e}lica Torres and Elisenda Feliu, \emph{Symbolic proof of bistability in
  reaction networks}, SIAM Journal on Applied Dynamical Systems \textbf{20}
  (2021), no.~1, 1--37.

\bibitem[TNP18]{Toth}
J{\'a}nos T{\'o}th, Attila~L{\'a}szl{\'o} Nagy, and D{\'a}vid Papp,
  \emph{Reaction kinetics: exercises, programs and theorems}, Springer, 2018.

\bibitem[Vas23]{Vas}
Nicola Vassena, \emph{Symbolic hunt of instabilities and bifurcations in
  reaction networks}, Discrete and Continuous Dynamical Systems-B (2023), 0--0.

\bibitem[Vas24]{Vas24noHurw}
\bysame, \emph{Mass action systems: two criteria for {H}opf bifurcation without
  hurwitz}, arXiv preprint arXiv:2402.18188 (2024).

\bibitem[VdDW02]{Van}
Pauline Van~den Driessche and James Watmough, \emph{Reproduction numbers and
  sub-threshold endemic equilibria for compartmental models of disease
  transmission}, Mathematical biosciences \textbf{180} (2002), no.~1-2, 29--48.

\bibitem[VdDW08]{Van08}
P~Van~den Driessche and James Watmough, \emph{Further notes on the basic
  reproduction number}, Mathematical epidemiology, Springer, 2008,
  pp.~159--178.

\bibitem[VG16]{Vyska}
Martin Vyska and Christopher Gilligan, \emph{Complex dynamical behaviour in an
  epidemic model with control}, Bulletin of mathematical biology \textbf{78}
  (2016), no.~11, 2212--2227.

\bibitem[Vol26]{volterra26}
Vito Volterra, \emph{Fluctuations in the abundance of a species considered
  mathematically}, Nature \textbf{118} (1926), no.~2972, 558--560.

\bibitem[VS24]{VasStad}
Nicola Vassena and Peter~F Stadler, \emph{Unstable cores are the source of
  instability in chemical reaction networks}, Proceedings of the Royal Society
  A \textbf{480} (2024), no.~2285, 20230694.

\bibitem[Wan06]{Wang}
Wendi Wang, \emph{Backward bifurcation of an epidemic model with treatment},
  Mathematical biosciences \textbf{201} (2006), no.~1-2, 58--71.

\bibitem[WG64]{MA64}
P~Waage and CM~Guldberg, \emph{Studier over affiniteten}, Forhandlinger i
  Videnskabs-selskabet i Christiania \textbf{1} (1864), 35--45.

\bibitem[WR04]{WR}
Wendi Wang and Shigui Ruan, \emph{Bifurcations in an epidemic model with
  constant removal rate of the infectives}, Journal of Mathematical Analysis
  and Applications \textbf{291} (2004), no.~2, 775--793.

\bibitem[XR07]{Ruan07}
Dongmei Xiao and Shigui Ruan, \emph{Global analysis of an epidemic model with
  nonmonotone incidence rate}, Mathematical biosciences \textbf{208} (2007),
  no.~2, 419--429.

\bibitem[XWJZ21]{xu2021complex}
Yancong Xu, Lijun Wei, Xiaoyu Jiang, and Zirui Zhu, \emph{Complex dynamics of a
  {SIRS} epidemic model with the influence of hospital bed number}, Discrete \&
  Continuous Dynamical Systems-B \textbf{26} (2021), no.~12, 6229.

\bibitem[YC18]{Yu}
Polly~Y Yu and Gheorghe Craciun, \emph{Mathematical analysis of chemical
  reaction systems}, Israel Journal of Chemistry \textbf{58} (2018), no.~6-7,
  733--741.

\bibitem[Zas95]{zaslavsky1995stability}
Boris Zaslavsky, \emph{Stability of {H}opf bifurcation in quasi-linear models},
  Journal of mathematical analysis and applications \textbf{192} (1995), no.~1,
  294--311.

\bibitem[ZF12]{ZhouFan}
Linhua Zhou and Meng Fan, \emph{Dynamics of an {SIR} epidemic model with
  limited medical resources revisited}, Nonlinear Analysis: Real World
  Applications \textbf{13} (2012), no.~1, 312--324.

\bibitem[ZL08]{zhang2008backward}
Xu~Zhang and Xianning Liu, \emph{Backward bifurcation of an epidemic model with
  saturated treatment function}, Journal of mathematical analysis and
  applications \textbf{348} (2008), no.~1, 433--443.

\bibitem[ZL23]{Zhang23}
Yingying Zhang and Chentong Li, \emph{Bifurcation of an {SIRS} model with a
  modified nonlinear incidence rate}, Mathematics \textbf{11} (2023), no.~13,
  2916.

\bibitem[ZXL07]{ZhouXiao}
Yugui Zhou, Dongmei Xiao, and Yilong Li, \emph{Bifurcations of an epidemic
  model with non-monotonic incidence rate of saturated mass-action}, Chaos,
  Solitons \& Fractals \textbf{32} (2007), no.~5, 1903--1915.

\end{thebibliography}
\bibliographystyle{amsalpha}

\end{document}